\documentclass[dvips]{article}

\usepackage{graphicx}

 \usepackage{amscd}
                     \usepackage{amssymb}
 \newcommand{\be}{\begin{equation}}
       \newcommand{\ee}{\end{equation}}
       \newcommand{\ba}{\begin{eqnarray}}
        \newcommand{\ea}{\end{eqnarray}}
 \newcommand{\ban}{\begin{eqnarray*}}
 \newcommand{\ean}{\end{eqnarray*}}

\newcommand{\vare}{\varepsilon}
 
 \newcommand{\lp}{\langle}
 \newcommand{\rp}{\rangle}
 \newcommand{\ra}{\rightarrow}

\newcommand{\R}{\Bbb{R}}
  \newcommand{\qed}{\hspace*{\fill}\rule{3mm}{3mm}\quad \vspace{.2cm}}
  \newcommand{\Pf}{\noindent {\bf Proof:} }

 \newcommand{\sect}[1]{\section{#1} \setcounter{equation}{0}}

 \newtheorem{theo}{Theorem}[section]
 \newtheorem{example}{Example}[section]

\begin{document}
 \newtheorem{defn}[theo]{Definition}
 \newtheorem{ques}[theo]{Question}
 \newtheorem{lem}[theo]{Lemma}
 \newtheorem{lemma}[theo]{Lemma}
 \newtheorem{prop}[theo]{Proposition}
 \newtheorem{coro}[theo]{Corollary}
 \newtheorem{ex}[theo]{Example}
 \newtheorem{note}[theo]{Note}
 \newtheorem{conj}[theo]{Conjecture}
 \newtheorem{rmk}[theo]{Remark}
 \newtheorem{rmrk}[theo]{Remark}

\title{Various Covering Spectra for Complete Metric Spaces}
\author{Christina Sormani \thanks{Partially supported by NSF Grant \# DMS-1006059 and a PSC CUNY Research Award
}
\and Guofang Wei \thanks{Partially
 supported by NSF Grant \# DMS-1105536.}}
 \date{}
 \maketitle

 \begin{abstract}
 Here we study various covering spectra for complete noncompact length
 spaces with universal covers (including Riemannian manifolds and the
 pointed Gromov Hausdorff limits of Riemannian manifolds with lower
 bounds on their Ricci curvature).   We relate the covering spectrum
 to the (marked) shift spectrum of such a space.   We define the slipping 
 group generated by elements of the fundamental group whose 
 translative lengths are $0$.   We introduce a rescaled length, the
 rescaled covering spectrum and the rescaled slipping group.   Applying
 these notions we prove that certain complete noncompact Riemannian
 manifolds with nonnegative or positive Ricci curvature
 have finite fundamental groups.   Throughout we suggest further
 problems both for those interested in Riemannian geometry and 
 those interested in metric space theory.   
 \end{abstract}

 \newcommand{\inj}{\mbox{inj}}
 \newcommand{\vol}{\mbox{vol}}
 \newcommand{\diam}{\mbox{diam}}
 \newcommand{\Ric}{\mbox{Ric}}
 \newcommand{\Iso}{\mbox{Iso}}
 \newcommand{\Hess}{\mbox{Hess}}
 \newcommand{\divg}{\mbox{div}}

\section{Introduction}

The covering spectrum of a compact Riemannian manifold or
length space captures the metric properties 
needed to obtain topological
information about the given space.    In prior work of the authors
\cite{SoWei3}, we applied the covering spectrum to determine
the properties of the fundamental group of a compact metric
space and to determine whether that space has a universal
cover.  We proved the covering spectrum is determined by the
marked length spectrum and the elements of the covering spectrum
lie in the half length spectrum.   Then de Smit, Gornet and Sutton
developed a means of producing pairs of compact manifolds which
have the same covering spectrum in \cite{deSmit-Gornet-Sutton}.
They produced pairs of compact manifolds with the same Laplace
spectrum that have different covering spectra in
\cite{deSmit-Gornet-Sutton} \cite{deSmit-Gornet-Sutton-2}.   
This is intriguing
in light of the work of Colin de Verdiere and Duistermaat-Guillemin
relating the length and Laplace spectra of compact Riemannian
manifolds \cite{CdV, DuGu}.   A recent extension of the notion of covering 
spectrum to a larger class of spaces which is called the critical spectrum has
been studied by Wilkins, Plaut, Conant, Curnutte, Jones, Pueschel and
Walpole \cite{Wilkins-thesis} \cite{Plaut-Wilkins:1205:1137} \cite{Wilkins-1208.3494} \cite{CCJPPWW}.  The key definitions,
theorems and examples are reviewed in Section 2.

In this paper we are concerned with the covering spectra and other
related spectra on complete noncompact Riemannian manifolds and
length spaces.   Prior work in this direction has been conducted by
the authors in \cite{SoWei4} where we developed the notion of
the cut off covering spectrum.  The cut off covering spectrum effectively
removed information about the space outside of large balls thus enabling 
us to extend a number of our prior results.   However data about the
behavior of the space at infinity was lost in the process.   

Here we
develop new spectra designed to capture the properties of 
complete noncompact Riemannian manifolds
at infinity.   More generally, we assume all of our
metric spaces, $X$, are complete noncompact length spaces
as in Definition~\ref{defn-complete-noncompact} and that they
have universal covers, $\tilde{X}$, in the sense that
the universal cover is a cover of all covering spaces (c.f. \cite{Sp}).
We do not assume the universal covering is simply connected.
Recall that the Gromov-Hausdorff limits of complete
noncompact Riemannian manifolds with uniform lower bounds on
their Ricci curvature were proven to have universal
covering spaces in \cite{SoWei2}.  So all results in this paper
apply to these essential limit spaces that have been explored
extensively by Cheeger, Colding, Ding, Ennis, Honda, Menguy,
Munn, Naber, Ohta, Tian, the authors and many others (c.f. 
\cite{ChCo1}\cite{CoNb}\cite{Ding}\cite{EnWei}\cite{Honda}\cite{Munn}\cite{Menguy}\cite{Ohta} ).

In the final 
section of the paper, we apply our spectra to prove special
cases of Milnor's conjecture that the fundamental group of
a complete noncompact manifold with nonnegative Ricci curvature
is finitely generated \cite{Mi}.  
Prior work in this direction
has been conducted by Li, Anderson, Wilking and the first
author \cite{Anderson}\cite{Li1986}\cite{Sormani-group}\cite{Wilking2000}. 
In fact we prove complete noncompact manifolds with certain
spectral properties that have
positive or nonnegative Ricci curvature have a finite fundamental group 
[Theorems~\ref{finite} and~\ref{Ricci-pi-rs-thm} respectively].
So one might prefer to view these as extensions of Myers' Theorem
\cite{Myers}
that compact manifolds with positive Ricci curvature have
finite fundamental groups.
Before we specialize to Riemannian manifolds, we first need
to extend our theorems about the covering spectrum proven
in \cite{SoWei3} to complete noncompact length spaces and introduce
our new covering spectra.

In Section 3, we extend our results in \cite{SoWei3}
relating the covering spectrum and
the length spectrum of a compact length space to the complete
noncompact setting.   We introduce the {\em (marked) shift spectrum} 
[Definition~\ref{defnshift}] which captures the lengths of elements
of the fundamental group in the complete noncompact setting and agrees
with the classical (marked) length spectrum on compact length spaces.   We 
prove Theorem~\ref{shiftdetcov}, that the the marked shift spectrum
determines the covering spectrum on such spaces (extending Theorem 4.7
of \cite{SoWei3}).   We prove the covering spectrum is a subset of the
closure of the half shift spectrum in Theorem~\ref{covofshift}.   In the
compact setting this was proven in \cite{SoWei3} without requiring a
closure, but in the complete noncompact setting we show this is necessary
with Example~\ref{excovofshift}.   We suggest
notion of a (marked) shift spectrum that might be studied on spaces without
universal covers [Remark~\ref{supinf}].   We close this section by
discussing elements of the fundamental group whose lengths are achieved
and proving Theorem~\ref{thm-line} that if infinitely many such elements
have their lengths achieved in a common compact set then the universal cover contains a line.  

In Section 4, we extend our results in \cite{SoWei3}
relating the covering spectrum
to the universal cover of a compact length space to the complete
noncompact setting.   We introduce the {\em universal slipping group},
$\pi_{slip}(X)$, and the {\em universal delta cover}, $\tilde{X}^0$
[Definitions~\ref{def-univ-slip} and~\ref{def-univ-delta-cov}].
On a compact metric space, the universal slipping group is trivial 
and the universal delta cover is the universal cover.  In
fact the authors proved that if the covering spectrum of a compact
length space has a positive infimum, then the universal cover 
exists and is a $\delta$ cover  \cite{SoWei3}.   
Here we present complete noncompact
manifolds whose $\delta$ covers are all trivial and whose
universal cover is nontrivial [Example~\ref{cuspcyl}].   Every
element of  the fundamental group in this example 
is represented by a sequence of loops
which slip out to infinity as their lengths decrease to $0$. 
We prove the universal delta cover, $\tilde{X}^0$, is a pointed 
Gromov-Hausdorff limit of $\delta$ covers
as $\delta \to 0$ [Theorem~\ref{GHcover}].   It need not
be a delta cover itself [Remark~\ref{rmrk-univ-del-cov-not-del-cov}],
nor is it the universal cover [Remark~\ref{rmrk-univ-del-cov-not-univ-cov}].
However, if the covering spectrum has a positive infimum,
then $\tilde{X}^0$ is a delta cover [Proposition~\ref{univ-del-cov-inf}].
In general we show $\tilde{X}^0=\tilde{X}/\pi_{slip}(X)$
[Theorem~\ref{univcov}].  Thus if the slipping group is empty
and the covering spectrum has a positive infimum, the
universal cover is a $\delta$ cover [Corollary~\ref{coro-univ-cov}].   
In Remark~\ref{open-univ-cov} we suggest that an adaption of
this work to complete length spaces which are not known to 
have universal covers might be applied to prove the existence of 
a universal cover for such a space.

In Section 5, we introduce two new scale invariant covering spectrum:
the basepoint dependent {\em rescaled covering spectrum}, 
$CovSpec_{rs}^{x_0}(X)$,
and the {\em infinite rescaled covering spectrum}, $CovSpec^\infty_{rs}(X)$
[Definitions~\ref{rescaledcovspec}].
Both these spectra are defined only on complete noncompact spaces
using $\delta$ rescaled covering groups [Definition~\ref{defn-rescaled-delta-group}] whose elements have rescaled lengths [Definitions~\ref{defn-rescaled-length} and~\ref{defn-infinite-rescaled-length}] defined by
measuring the lengths sequences of representative loops based at
points diverging to infinity.   These rescaled covering spectra 
take values lying in $[0,1]$.
We also define the {\em rescaled slipping group}, $\pi_{rs}^\infty(X,0)$,
which is generated by elements of the fundamental group that have representative loops diverging to infinity whose rescaled lengths 
converge to $0$. In the final subsection we suggest a few
means of defining these spectra on spaces without universal covers.  Nevertheless, we believe the rescaled covering spectra will prove most useful in the study of submanifolds of Euclidean space and Riemannian manifolds, both of which always have universal covers.

In Section 6, we study the asymptotic behavior of complete noncompact
metric spaces using the rescaled covering spectra and the
rescaled slipping group.  We prove that if the rescaled length of 
an element of the fundamental group is strictly less than $2$, then 
it is represented by a sequence of loops diverging to infinity [Lemma~\ref{rsloopstoinfty}].    We prove that
if  $CovSpec_{rs}^\infty(X)\subset [0,1)$
then the cut off covering spectrum of \cite{SoWei4}
is trivial [Theorem~\ref{rescalevscut}].  These
two spectra are recording very different metric topological information
about a space.  Next we compute the rescaled covering spectra of
one sheeted hyperboloids [Example~\ref{one-sheet}]
and doubly warped products [Theorem~\ref{asymthm}].   We close with 
a conjecture about the rescaled covering spectra of manifolds with linear
diameter growth [Conjecture~\ref{conj-diam-growth}].
 
The notions and theorems in Sections 2-6 are then applied in the
final section to complete noncompact Riemannian manifolds with
curvature bounds.   First we study manifolds with nonnegative
sectional curvature, whose topology can be understood because such
manifolds have a compact soul.   We prove a noncompact
versions of Myers' Theorem: 

\begin{theo} \label{finite}
If $M$ has positive Ricci curvature
and $K\subset M$ is compact,
then only finitely many elements of the fundamental group
have their lengths achieved within $K$.
\end{theo} 

This theorem is proven by constructing a line in the universal cover
and applying the Cheeger-Gromoll Splitting Theorem \cite{ChGr1971}.   
Note
that in the classic example of Nabonnand
\cite{Nab}, which has positive Ricci curvature,
there are infinitely many elements in the fundamental
group as well, but their lengths are not achieved anywhere.   
In a cylinder, which has nonnegative Ricci curvature,
infinitely many elements can have their length achieved
in a common compact set however the rescaled covering spectrum
of the cylinder is trivial and the fundamental group lies
in the rescaled slipping group.    In fact we prove:

\begin{theo} \label{Ricci-pi-rs-thm}
If $M^n$ is a complete noncompact Riemannian manifold whose
rescaled slipping group, $\pi_{rs}^\infty(M,0)$, is trivial and whose rescaled
covering spectrum, $CovSpec_{rs}^\infty(M)$, has a positive infimum, then
it has a finite fundamental group.
\end{theo}

This theorem is proven by applying the Bishop-Gromov Volume
Comparison Theorem in a style similar
to Milnor's proof that finitely generated subgroups of the fundamental 
group have polynomial growth \cite{Mi}.   Note that in the example of
Nabonnand the entire fundamental group lies in the rescaled
slipping group.
We observe that a similar statement holds when $M$ is the pointed
Gromov-Hausdorff limit of manifolds with nonnegative Ricci curvature
[Remark~\ref{lim-Ricci-pi-rs-thm}] including the tangent cone at
infinity for such a space.   

It is possible that the techniques applied to prove these theorems may
lead to a proof of the Milnor Conjecture that the fundamental group
of a manifold with nonnegative Ricci curvature is finitely generated.
Our proofs apply the Cheeger-Gromoll Splitting Theorem and a
Milnor style application of the
Bishop Gromov Volume Comparison Theorem.   Other techniques for
controlling the fundamental groups of such manifolds that appear 
in work of Anderson, Li, Wilking and the first author in \cite{Anderson}\cite{Li1986}\cite{Sormani-group}\cite{Wilking2000}
have not yet been applied in combination with our new ideas.
Further suggested problems are stated in the final remarks of the
paper. 

The authors would like to thank Ruth Gornet, Carolyn Gordon 
and Shing-Tung Yau
for their interest in the covering spectrum and Burkhard Wilking for 
suggesting an analysis of ${\Bbb{S}}^3\times\Bbb{R}^4/ \textrm{Pin}(2)$.

\section{Background}

In this section we review our prior work and that of others
taking advantage of the simplified definitions one can use when there 
is a universal cover.

 \begin{defn} \label{defn-complete-noncompact}
{\em
A {\em complete length space} is a
complete metric space such that every pair of points in the space is joined by
 a length minimizing rectifiable curve.  The distance between
 the points is the length of that curve.  A {\em compact length
 space} is a compact complete length space.  (c.f. \cite{BBI}).
}
 \end{defn}

Note that complete Riemannian manifolds are complete length spaces
by the Hopf-Rinow Theorem.   As in Riemannian geometry, we
define geodesics as follows (c.f. \cite{BBI} \cite{Sormani-length}):

\begin{defn}
A geodesic $\gamma: I \to X$, is a locally length minimizing curve: 
\be
\forall \,t \in I \,\,\,\exists \,\vare_t>0 \,\,\,such \,\,\, that\, \,\, 
d_X\big(\gamma(t-\vare), \gamma(t+\vare)\big)= L\big(C([t-\vare, t+\vare])\big).
\ee
\end{defn}

The length spectrum of a Riemannian manifold in a classical
notion studied by many people over the years (c.f. \cite{CdV}).  
Here we recall the definition on a complete
length space (c.f.\cite{Sormani-length}):

\begin{defn}
The length spectrum is defined:
\be
Length(X)=\{L(\gamma): \,\, \gamma: {\Bbb{S}}^1 \to X\}
\ee
where $\gamma$ is a closed geodesic.
\end{defn}

Recall that on a metric space, $X$, with a universal cover, $\tilde{X}$, 
every element, $g$, in the fundamental group, $\pi_1(X)$, can be thought 
of as a deck transform of the univeral cover $g: X\to X$.
Thus every element has a well defined length
\be \label{eqn-length}
L(g) =  \inf _{x\in X} d_{\tilde{X}}(g\tilde{x},\tilde{x}) \in [0,\infty),  \label{L(g)}
\ee
where $\tilde{x}$ is an arbitrary lift of $x$ to the universal
cover.  
On complete noncompact 
spaces the length in (\ref{eqn-length}) might not be
achieved and could be $0$.  It is sometimes referred to as the
{\em translation length} in $\textrm{CAT}(0)$ geometry
(c.f. \cite{BO}).
On compact manifolds this infimum is achieved at some $x_0$ and
any minimal geodesic running between $\tilde{x}_0\in X$ and $g\tilde{x}_0$
can be seen to project to a geodesic loop, $\sigma:{\Bbb{S}}^1 \to X$.

\begin{defn}
The  marked length spectrum is the map:
\be
L: \pi_1(X) \to [0,\infty)
\ee
where $L(g)$ is defined as in (\ref{eqn-length}).   The
image $L(\pi_1(X))$ is a subset
of the length spectrum on a compact manifold because it
is achieved by a closed geodesic.  
\end{defn}

Let $\pi_1(X,\delta)$ be the subgroup of the
fundamental group generated by elements of length $< 2\delta$:
\be \label{defn-pi-delta}
\pi_1(X,\delta)= \lp g \in \pi_1(X): L(g) <2\delta \rp. 
\ee
In \cite{SoWei4} Thm 2.14, the
authors prove that the $\delta$ cover, $\tilde{X}^\delta$ of a 
metric space with
a universal cover can be found by taking
\be \label{groupcover} \label{defn-delta-cover}
\tilde{X}^\delta= \tilde{X}/ \pi_1(X, \delta).
\ee
 
 The covering spectrum of metric space (initially defined
 in \cite{SoWei3}) measures the size of
 one dimensional holes in the space: 
  
\begin{defn}  \label{def-cov-spec}
 Given a complete length space $X$, the covering spectrum of
 $X$, denoted  CovSpec$(X)$ is the set of all $\delta > 0$ such
 that
 \be \label{covspecdef1}
 \tilde{X}^\delta \neq \tilde{X}^{\delta'}
 \ee
 for all $\delta'> \delta$.
   \end{defn}
   
When $X$ has a universal cover then $\delta \in CovSpec(X)$
iff
 \be \label{covspecdef2}
 \pi_1(X, \delta) \neq \pi_1(X, \delta')
 \ee
 for all $\delta'> \delta$.   Note that $\delta$
covers are monotone in the sense that if $\delta_1<\delta_2$ then
$\tilde{X}^{\delta_1}$ covers $\tilde{X}^{\delta_1}$.  In a
compact manifold, where all lengths of elements are positive, 
any covering space is covered by a $\delta$
cover \cite{SoWei1}.   In \cite{SoWei3},
the authors proved that the covering spectrum of a 
compact length space is a subset of the half length
spectrum:  
\be\label{cov-to-length}
\lambda \in CovSpec(X) \qquad \Longrightarrow \qquad
2\lambda \in Length(X).   
\ee

\begin{example} \label{torus}
Let $X$ be the $k$ torus created by taking the isometric product of
$k$ circles:
\be
X={\Bbb{S}}^1_{r_1}\times {\Bbb{S}}^1_{r_2} \times \cdots \times {\Bbb{S}}^1_{r_k}
\ee
where $r_1\ge r_2\ge\cdots \ge r_k$ and ${\Bbb{S}}^1_r$ denotes a circle
of intrinsic diameter $r$ (circumference $2r$).
Then if $\delta \in (r_{j-1}, r_{j}]$ we have
\be
\tilde{X}^\delta= \mathbb{R}\times\mathbb{R} \times \cdots
\times \mathbb{R}\times {\Bbb{S}}^1_{r_j} \times \cdots \times {\Bbb{S}}^1_{r_k}.
\ee
The covering spectrum 
is $\{r_1,r_2,r_3, \cdots, r_k\}$. 
\end{example}

\begin{example}\label{fathandle}
If $X$ is a genus 2 surface which is a very fat 
figure eight with two small holes, the covering spectrum
has two elements $\{\lambda_1,\lambda_2\}$ where $\lambda_1$ is 
half the length of the shortest geodesic around the smaller hole
and $\lambda_2$ is 
half the length of the shortest geodesic around the other hole.
These two elements then generate the entire fundamental group.
\end{example}

Bart de Smit, Ruth Gornet and Craig Sutton have developed a
method of producing two compact manifolds with the same covering
spectrum, building upon work of Sunada used to produce
pairs of compact manifolds with the same Laplace spectrum in 
\cite{deSmit-Gornet-Sutton}.   Despite the close
relationship between the Laplace and Length spectrum proven by
Colin de Verdiere in \cite{CdV} and the close relationship
between the Covering and Length spectrum proven by the authors
in \cite{SoWei2}, de Smit, Gornet and Sutton have found 
pairs of compact Laplace isospectral  manifolds with different covering spectra.   They study higher dimensional manifolds in \cite{deSmit-Gornet-Sutton} and surfaces in \cite{deSmit-Gornet-Sutton-2}.

In the complete noncompact setting, where the lengths of
elements of the fundamental group are not achieved by lengths
of closed geodesic loops, one still seems to have some relationship
between the covering spectrum and the translative length spectrum,
as indicated by the following example:

\begin{example} \label{ex-loop-to-one}
Let $M^2$ be the warped product $\Bbb{R}\times_f {\Bbb{S}}^1$
where $f(x)=1 + e^{-x^2}$.  For $\delta\le \pi$, all
balls of radius $\delta$ are simply connected
because $f(x)\diam({\Bbb{S}}^1)\ge \pi$ for all $x$.   Thus the
$\delta$ cover is the universal cover,
$\tilde{M}=\Bbb{R}\times_f \Bbb{R}$, for
$\delta \le \pi$.  For $\delta>\pi$,
we capture a loop $\beta$ running around ${\Bbb{S}}^1$ in a ball
of radius $\delta$ where $x\circ \beta =\sqrt{-\log(\delta/\pi-1)}$.  
So the $\tilde{M}^\delta=M$ for $\delta>\pi$.  Thus the 
covering spectrum is just $\{\pi\}$.  On the other hand one can show
there is no closed geodesic of length $2\pi$ in the
marked length spectrum, because any path $\gamma$
which traverses around the ${\Bbb{S}}^1$, has length
\be
L(\gamma)\ge f(\max \{x\circ\gamma\}) 2\pi>2\pi.
\ee
Nevertheless $L(g)=2\pi$ where $g$ is the deck transform
generating the fundamental group.
\end{example}

In the next section we prove Theorems~\ref{shiftdetcov}
and Theorem~\ref{covofshift},
relating the lengths of elements of the fundamental
group to the covering spectrum for all complete length
spaces.

\section{The Covering Spectrum and the Shift Spectrum}

In this section we discuss the relationship between the
covering spectrum and the length spectrum.
We make the following new definition which agrees with the
marked length spectrum on a compact length space.

\begin{defn} \label{defnshift}
The shift spectrum, $Shift(X)$,
is the collection of translative
lengths $L(g)$ of elements of the fundamental group $\pi_1(X)$.
The marked shift spectrum is the map which
takes each element of the fundamental group $\pi_1(X)$
to its length.
Two spaces $X$ and $Y$ are said to have the same marked
shift spectrum if there is an isometry between their
fundamental groups, $f:\pi_1(X)\to\pi_1(Y)$
such that $L(f(g))=L(g)$ for all $g\in \pi_1(X)$.
\end{defn}

It has been suggested that this should be called the
translative length spectrum.  Since ``shift'' is a shorter
word than ``translative'' and has the same meaning, we use
it instead.

As seen in Example~\ref{ex-loop-to-one}, the shift spectrum is
not necessarily a subset of the length spectrum on a complete
noncompact manifold.  In that example, $2\pi \in Shift(M)$,
but is not the length of any smoothly closed geodesic.

\subsection{The marked shift spectrum determines the covering spectrum}

We can now extend Theorem 4.7
of \cite{SoWei3} to complete noncompact spaces.

\begin{theo} \label{shiftdetcov}
The marked shift spectrum of a complete
length space, $X$, with a universal cover
determines the covering spectrum of $X$.
\end{theo}

\Pf From Definition~\ref{def-cov-spec},  for $X$ with a universal cover,
 the covering spectrum is given by 
the $\delta$ where the covering groups $\pi_1(X,\delta)$
change, where is $\pi_1(X,\delta)$ is the subgroup of the
group of deck transforms generated by the deck transforms
of length $<2\delta$.

So if two spaces, $X$ and $Y$, share the same shift spectrum,
\be\label{Gnaut}
\pi_1(X,\delta) =\pi_1(Y,\delta)
\ee
via the restricted isometry between $\pi_1(X)$ and $\pi_1(Y)$.

Suppose on the contrary
that $X$ and $Y$ have the same marked shift spectrum and
different covering spectra.
We may assume without loss of generality that
that $\delta_0 \in CovSpec(X) \setminus CovSpec(Y)$.
Since
$\delta_0 \notin CovSpec(Y)$, then by definition, there is some
$\delta_1>\delta_0$ such that
\be
 \tilde{Y}^{\delta_0} = \tilde{Y}^{\delta_1}.
 \ee
Thus $\pi_1(Y,\delta_1)=\pi_1(Y,\delta_0)$.  Applying (\ref{Gnaut})
we get
\be
\pi_1(X, \delta_1)=\pi_1(Y,\delta_1)=\pi_1(Y,\delta_0)=\pi_1(X, \delta_0).
\ee
Thus we have
\be
 \tilde{X}^{\delta_0} = \tilde{X}^{\delta_1}.
 \ee
and so $\delta_0 \notin CovSpec(X)$ which is a contradiction.
\qed

\subsection{The covering spectrum lies in the closure of the shift spectrum}

Recall the lower semiclosure of a set $A\subset \Bbb{R}$, denoted
$Cl_{lower}(A)$, is the set of all limits of nonincreasing
sequences of points in $A$ (c.f. \cite{SoWei4}).

\begin{theo}\label{covofshift}
The covering spectrum is a subset of the lower semiclosure of the 1/2
shift spectrum:
\be
CovSpec(X) \subset Cl_{lower}(\{h/2: h \subset \mbox{Shift} \, (X)\}).
\ee
\end{theo}

Before proving this, we provide an example demonstrating that
the covering spectrum is not a subset of the half shift spectrum:

\begin{example}\label{excovofshift}
Let $X$ be a collection of circles of intrinsic
diameter $\{\pi+\pi/j: j\in \Bbb{N}\}$
joined at a common point.  Let $g_j \in \pi_1(X)$ be
the element represented by a loop going once around the
circle of intrinsic diameter $\pi+\pi/j$.   Then
\be
\pi_1(X, \delta)=\lp g_k, g_{k+1}, g_{k+2}, ... \,\, : \pi+\pi/k<\delta \rp
\ee 
and
\be
CovSpec(X)=\{\pi(1+1/j): j\in \Bbb{N}\} \cup \{\pi\}.
\ee
Here we have included $\pi$ because for $\delta= \pi$,
$\pi_1(X, \delta)$ is trivial, and for $\delta'>\pi$, there
exists $k$ sufficiently large that $\pi+\pi/k<\delta'$ and
$\pi_1(X, \delta')$ is nontrivial.
On the other hand,
the shift spectrum is the collection of finite sums:
\be
Shift(X)=\left\{\sum_{j=1}^N k_j 2\pi(1+1/j): \, k_j \in \{0\} \cup \Bbb{N}, \,
N \in \Bbb{N} \right\}
\ee
which contains twice every value in $CovSpec(X)$ except 
the value $\pi$.
\end{example}

Keeping this example in mind, we now prove Theorem~\ref{covofshift}:
\vspace{.3cm}

\Pf
Suppose on the contrary that there is a space $X$ and an
element
\be
\delta_0\in CovSpec(X) \setminus Cl_{lower}((1/2)Shift(X)).
\ee
Then there exists $\epsilon>0$ such that
\be
[\delta_0+\epsilon) \cap (1/2)\mbox{Shift}\,(X) =\emptyset.
\ee
By Definition~\ref{defnshift},
\be
\forall g \in \pi_1(M), L(g)\neq [2\delta_0, 2\delta_0+2\epsilon).
\ee
Thus,
\be
\pi_1(X,\delta_0+\epsilon/2)=\pi_1(X,\delta_0), 
\ee
and
 $\delta_0$ is not in the covering spectrum.
\qed

\begin{coro}\label{covsubshift}
If $X$ has a closed shift spectrum, then
\be
CovSpec(X) \subset \mbox{Shift}\, (X).
\ee
\end{coro}

\begin{rmk} \label{supinf}
It is possible one might achieve the same results for spaces which do
not have universal covers possibly by extending the definition in (\ref{L(g)})
as follows:
\be
L(g) = \sup_{\tilde{X}} \inf_{C} d_{\tilde{X}}(\tilde{C}(0),\tilde{C}(1))
\ee
where the supremum is taken over all regular covering spaces, $\tilde{X}$ 
of $X$ and
the infimum is taken over all loops $C:[0,1]\to X$ freely homotopic to
a representative of $g$, and where $\tilde{C}$ is a lift of $C$ to
$\tilde{X}$.  This leads to a natural extension of the notion of the shift
spectrum.  Proving the extensions would be difficult without a common
cover to examine, however one might examine how we overcome this
issue in the proof that the covering spectrum is contained in the
half length spectrum in \cite{SoWei3} where no universal cover is assumed.
\end{rmk}

\subsection{Lines in universal covers}

Naturally, in a complete noncompact Riemannian manifold, some elements of
the shift spectrum may lie in the length spectrum.    This occurs, for
example, when the lengths of elements of the fundamental group
are actually achieved within a compact set.   The following theorem
demonstrates that this is in some sense exceptional:

\begin{theo}\label{thm-line}
Given a complete Riemannian manifold $M$, if there are infinitely many distinct 
elements $g_1, g_2, ..., g_i, ... $ of $\pi_1(M)$ such 
that all $L(g_i)$ are achieved in a compact set, then the universal 
cover of $M$ has a line.
\end{theo}

\Pf 
Let $\gamma_i$ be a representative of $g_i$ 
such that $L(\gamma_i)=L_i =l(g_i)$ and $\gamma_i(0)\in K$.  Since all $g_i$ are distinct, there is a subsequence $L_j \to \infty$. Fix a lift of $K$, $\tilde{K}$ in the universal cover. Let 
$\tilde{\gamma}_j$ be a lift of $\gamma_j$ with $\tilde{\gamma}_j(0) \in \tilde{K}$. Since $L(\gamma_j)=l(g_j)$, $g_j$ is a translation along $\tilde{\gamma}_j$, and   $\tilde{\gamma}_j$ is minimal on any subinterval of length $L_j$.  In particular $\tilde{\gamma}_j$ is minimal on $[-L_j/2, L_j/2]$. Since $\tilde{\gamma}_j(0) \in \tilde{K}$ and $\tilde{K}$ is compact, $\tilde{\gamma}_j$ converges to a line.
\qed

\section{The Covering Spectrum and the Universal Cover}

On a compact length space, $\delta$ covers were used to
prove the existence of a universal cover.  In fact the
universal cover is a $\delta$ cover for a sufficiently
small $\delta$ \cite{SoWei1}, and the universal
cover of a compact length space exists iff its covering
spectrum is finite \cite{SoWei3}.
The same is not true in the complete
noncompact setting, even when the space is a Riemannian
manifold:

\begin{example} \label{cuspcyl}
  The cylindrical cusp manifold
\be 
M=\Bbb{R} \times_{e^r} {\Bbb{S}}^1
\ee
has $\tilde{M}^\delta=M$ for all $\delta$ while the universal
cover is diffeomorphic to a plane.
\end{example}

In this section we investigate the relationship
between the delta covers and the universal covers of
a complete noncompact manifold.  Note that in \cite{SoWei2}
we used relative $\delta$ covers of balls to prove existence of
universal covers in some noncompact settings, but here
we are focusing on the delta covers of the manifold itself.

\subsection{The fundamental group and the universal slipping group}

We begin by examining the elements of the fundamental group
which cause the problem seen in Example~\ref{cuspcyl}.

\begin{defn} 
For each element $g \in \pi_1(M)$, $p \in M$, 
let $L(g,p) =d_{\tilde{M}} (\tilde{p}, g \tilde{p})$, 
where $\tilde{p}$ is some lift of $p$ in the 
universal cover $\tilde{M}$, i.e. the length of a 
shortest representative of $g$ at $p$.
\end{defn}

\begin{defn} \label{defn-slipping}
The slipping elements $g$ in the fundamental group
such that
\be
L(g)=\inf_{x\in \tilde{M}} d_{\tilde{M}}(x, gx)=0  
\ee
The loops representing these elements
slide out to infinity and their lengths disappear.
\end{defn}

The set of slipping elements is not a group since it need not be closed
as we will see in Example~\ref{ex-4.2}.
However, it generates a group we call the {\em slipping
group}.

\begin{example} \label{ex-4.1}
The slipping group of the isometric product cylinder $\Bbb{R} \times_1 {\Bbb{S}}^1$
is empty.  The slipping group of a warped product cylinder
$\Bbb{R} \times_f {\Bbb{S}}^1$ where $\lim_{r\to\infty}f(r)=0$ is the entire fundamental
group, ${\Bbb{Z}}$.  The slipping group of the doubly warped product
$\Bbb{R} \times_f {\Bbb{S}}^1 \times_h {\Bbb{S}}^1$ where $f$ is uniformly bounded below
but $\lim_{r\to\infty}h(r)=0$  is ${\Bbb{Z}}$ which is a subgroup of
$\pi_1={\Bbb{Z}}\times {\Bbb{Z}}$ .
\end{example}

Of course in these examples the slipping elements and the slipping
group agree.  Next we see why we need to generate the slipping group.

\begin{example}  \label{ex-4.2}
We begin by taking a figure 8 formed by taking two circle of radius
1 joined at a point.  Let $M^2$ be a smooth surface of genus 2 which
is obtained by taking the boundary of the tubular neighborhood of
radius $1/2$ about this figure eight in $\Bbb{R}^3$ and smoothing it slightly.
Cross $M^2$ with a line.
For the negative direction on the line, $r \in (-\infty,1]$,
take the isometric product metric.  For the positive direction,
$r\in [1,\infty)$,
of the line, change the metric on $M^2$ smoothly so that
it is the smoothened tubular neighborhood of radius $1- 1/(2r)$
about the fixed figure eight.

Now the fundamental group of the figure eight (and $M^2 \times \Bbb{R}$)
is the free group on two elements $g_1$ and $g_2$ where $g_1$ goes once
around the first circle and $g_2$ goes once around the second.  The slipping
elements of $M^2 \times \Bbb{R}$ with this metric, are $g_1^j$ and $g_2^j$
where $j\in {\Bbb{Z}}$.  The slipping group is the whole fundamental group.
However $g_1g_2$ and other mixed elements are not slipping elements, because
in order to go around both holes, they must have length $\ge 2$.
\end{example}

There is a possibly larger group we can call the
{\em universal slipping group} which is defined as follows:

\begin{defn} \label{def-univ-slip}
An element $g$ in the fundamental group is in the universal
slipping group of $M$, denoted $\pi_{slip}(M)$,
if for all $\delta>0$ there exist elements
$g_1, g_2, ...g_N \subset \pi_1(M)$ such that $L(g_j)<\delta$
and $g=g_1g_2\cdots g_N$.  Here $N$ may depend on $\delta$
and on $g$.
\end{defn}

Clearly the slipping group is a subgroup of the universal slipping
group.  The universal slipping group may be strictly larger than the
slipping group as can be seen in the following example:

\begin{example} \label{pants} 
We construct a complete noncompact surface, $M^2$, as follows.
First we take the ``pair of pants'', $P^2$, and endow it
with a Riemannian metric so that the region near the ``waist'' of the pants,
is isometric to a cylinder $S_{2\pi}^1\times [0,\vare)$ and the regions
near the legs are isometric to cylinders $S_{\pi}^1\times [0,\vare)$ where
${\Bbb{S}}^1_{r}$ is a circle of intrinsic diameter $r$.   We further require that the 
shortest loop homotopic to the waist in $P^2$ has length $4\pi$
and the shortest loop homotopic to a leg in $P^2$ has length $2\pi$.
We glue together two pairs
of pants at the waist and call the waist where they are glued $\gamma_1$
.  Then we glue on four more pairs of pants (rescaled
by $1/2$) so that their waists are glued into the four legs
along $\gamma_{2,1}, \gamma_{2,2}, \gamma_{2,3}, \gamma_{2,4}$.   
Next we glue
eight more pairs of pants (rescaled by $1/4$) so that their waists are glued to the eight legs along $\gamma_{3,1},...\gamma_{3,8}$
and so on ad infinitum.   This forms a complete noncompact
metric space with a collection of closed geodesics $\gamma_{i,j}$
where $i \in \mathbb{N}$ and $j\in 1, \cdots, 2^{i}$.   These geodesics are
the shortest curves in their free homotopy classes.   If we set
$\sigma_{i,j}$ to be a minimizing geodesic from $\gamma_1(0)$
to $\gamma_{i,j}(0)$, then we have 
\be
g_{i,j}=[ \sigma_{i,j}^{-1}* \gamma_{i,j}* \sigma_{i,j}] \in \pi_1(M^2)
\ee
such that $L(g_{i,j}) = L(\gamma_{i,j}) = (1/2)^i 2\pi$ and
(assuming they are oriented and ordered consistently)
\be
g_{i,j}=g_{i+1,2j-1} g_{i+1,2j}
\ee
Thus every element of $\pi_1(M^2)$ is in the universal slipping
group although no element of $\pi_1(M^2)$ is in the slipping group.
\end{example}

We have the following nice description for the universal
slipping group:  

\begin{lemma}\label{pideltaslip}
When viewed as subsets of $\pi_1(M)$:
\be
\pi_{slip}(M) = \bigcap_{\delta>0} \pi_1(\tilde{M}^\delta).
\ee
\end{lemma}

\Pf Let $g\in \bigcap_{\delta>0} \pi_1(\tilde{M}^{\delta})$, then for all $\delta>0$
$g\in \pi_1(\tilde{M}^{\delta})$.   By the definition of $\pi_1 (M, \delta)$
in (\ref{groupcover}), we know there
exists elements of length $<2\delta$ whose product is $g$, and since
this is for all $\delta$, g is in the slipping group.

Assume on the other hand that g is in the slipping group, then
for all $\delta >0, \ g$ is a product of elements of length less than $\delta$
so it is in $\pi_1(M, \delta/2)$.
\qed

\subsection{The universal cover and the universal $\delta$ cover}

On a compact length space with a universal cover we 
proved in \cite{SoWei3} that
the covering spectrum is finite and that the universal cover is
the $\delta$ cover corresponding to the smallest $\delta$ in the
spectrum.  On a complete noncompact length space with universal cover,
clearly the covering spectrum may have infinitely many elements:

\begin{example}  \label{surface-infinity-genus}
A surface of infinite genus with the same size holes has a finite
covering spectrum.  However, if the size of the holes goes
to zero as they approach to infinity, then the surface has a covering
spectrum consisting of infinitely many elements whose infimum is 0
and the universal cover is not a $\delta$ cover.
\end{example}

Note that in this example the universal cover is not a delta cover
but it is a limit of the $\delta$ covers as $\delta \to 0$.  That
is, $\tilde{M}$ is the pointed Gromov Hausdorff limit of
$\tilde{M}^\delta$ as $\delta \to 0$ where the points at the centers
of the balls used for this limit are the lifts of a fixed point in the
space.

We introduce the universal delta cover and prove that it can
always be obtained as such a pointed Gromov-Hausdorff
limit:

\begin{defn} \label{def-univ-delta-cov}
The universal delta cover, $\tilde{M}^0$, 
of a length space, $M$, is a covering space
which covers all delta covers, ${\tilde M}^\delta$.  We require further that 
for any other covering space $\tilde{M}'$ which covers all $\delta$
covers, ${\tilde M}^\delta$, we have $\tilde{M}'$ covers $\tilde{M}^0$.
\end{defn}

\begin{theo} \label{GHcover}
The pointed Gromov-Hausdorff limit,
\be
\lim_{\delta\to 0}\tilde{M}^\delta
\ee
of  a sequence of $\delta$ covers of a fixed complete
length space with universal cover, $M$, based at the lifts of
a fixed point $p\in M$ exists and is the universal $\delta$ cover
$\tilde{M}^0$.  Thus it doesn't depend on the basepoint used
to define the pointed Gromov-Hausdorff limit. 
\end{theo}

\Pf
First we will just show that for any $p$ and any sequence
$\delta_j \to 0$, there is a converging subsequence which converges
to a cover.  To get uniqueness, we will show it is the
universal delta cover in the sense described.

Let $f_j: \tilde{M}^{\delta_j} \to M$ be the covering maps so that
$f(p_j)=p$.  Let $h_j: \tilde{M} \to \tilde{M}^{\delta_j}$ be
covering maps so that $f(\tilde{p})=p_j$.
For any fixed $R, \epsilon>0$  let $N(\epsilon,R)$ be the maximum number of
disjoint balls of radius $\epsilon$ in $B_{\tilde{p}}(R) \subset \tilde{M}$.
This provides a uniform upper bound on the number of disjoint
balls of radius $\epsilon$ in $B_{p_j}(R) \subset \tilde{M}^{\delta_j}$.
Thus by Gromov's compactness theorem a subsequence of the $B_{p_j}(R)$
converges.  By Grove Petersen's Arzela Ascoli theorem and the fact that
$f_j$ and $h_j$ are distance nonincreasing maps, subsequences of these
functions converge as well.  Taking $R \to \infty$ and diagonalizing
we get a limit space $\tilde{M}^0$ with distance nonincreasing maps
$h_\infty:\tilde{M} \to \tilde{M}^0$ and $f_\infty: \tilde{M}^0 \to M$
such that the concatenation of these functions is a covering map.
Thus $f_\infty$ must be a covering map as well since it will act as
an isometry on any ball lifted isometrically to the universal cover.

Now this covering space we have obtained (which may depend on the
subsequence $\delta_j$ and on the point $p$), covers all the $delta$ covers
of the space.  This can be seen, by fixing $\delta>0$ and
observing that eventually $\delta_j<\delta$, so it is a limit of
spaces which cover $\tilde{M}^\delta$, and thus the same arguement
used above to explain why it covers $M$ can be used to explain why it
covers $\tilde{M}^{\delta_j}$.

Furthemore if $N$ covers all of $\tilde{M}^\delta$, then $N$ covers
all the $\tilde{M}^{\delta_j}$ and, as argued above these covering
maps have a convergence subsequence to a covering map from $N$
to $\tilde{M}^0$.

Thus $\tilde{M}^0$ is the universal $\delta$ cover
and so it is unique and doesn't depend on the choice of $p$
or the subsequence, and so no subsequence was required for the
limit after all.
\qed

We can now extend Theorem 3.4 of \cite{SoWei3}: 

\begin{prop} \label{univ-del-cov-inf}
If the infimum of the covering spectrum is positive
then the universal delta cover is a delta cover and the
infimum is in the covering spectrum.
\end{prop}

\Pf
If the spectrum has a positive infimum $\inf \{\delta\}$,
let $\delta_0=\inf \{\delta\}/2$.  Then
the  $\delta_0$ cover will cover all $\delta$ covers and so it must be the
universal delta cover by the uniqueness in Theorem~\ref{GHcover}.
Now there exists $\delta_i$ in the covering spectrum decreasing
to $\inf \{\delta\}$, and so by Defn~\ref{def-cov-spec}, we have
\be
\tilde{X}^{\delta_i} \neq \tilde{X}^{\delta'} \qquad \forall \delta'>\delta_i.
 \ee
By the monotonicity of $\delta$ covers, the $\inf \{\delta\}$ cover
must cover these $\delta_i$ covers and so for every $i$ we have
\be
\tilde{X}^{\inf\{\delta\}} \neq \tilde{X}^{\delta'} \qquad \forall \delta'>\delta_i.
 \ee
Taking $\delta_i$ down to $\inf\{\delta\}$ we see that
$\inf\{\delta\}$ satisfies the requirements of Defn~\ref{def-cov-spec}.
\qed

\begin{rmrk} \label{rmrk-univ-del-cov-not-del-cov}
In general, however, the universal delta cover need not be a delta cover.
This can be seen in Example~\ref{surface-infinity-genus}, where
each $\delta$ cover unravels finitely many holes and the
universal delta cover unwraps all of them.
\end{rmrk}

\begin{rmrk} \label{rmrk-univ-del-cov-not-univ-cov}
The universal delta cover
need not be the universal cover either.  In the
cusp cylinder of Example~\ref{cuspcyl}, 
$M$, all the delta covers, $\tilde{M}^\delta$,
are just isometric to $M$ and thus so is their Gromov-Hausdorff
limit, the universal delta cover.   Observe that in this example the
universal slipping group $\pi_{slip}$, is the whole fundamental group.
\end{rmrk}

\begin{theo} \label{univcov}
Given a complete length space, $M$, with a universal cover, $\tilde{M}$, fundamental group
$\pi_1(M)$, and universal slipping group $\pi_{slip}(M) \subset \pi_1(M)$.
Then the universal delta covering space $\tilde{M}^0$
satisfies
\be
\tilde{M}^0 =\tilde{M} / \pi_{slip}(M)
\ee
and is thus a regular covering space.
\end{theo}

\Pf  From  Lemma~\ref{pideltaslip}, it is enough to show \[ 
\tilde{M}^0 =\tilde{M} / (\bigcap_{ \delta >0}  \pi_1(M, \delta)).\]
This follows immediately since  $\tilde{M} / (\bigcap_{ \delta >0}  \pi_1(M, \delta))$ covers all $\delta$-cover $\tilde{M}^\delta$ and any cover covers all the $\delta$-cover covers $\tilde{M} / (\bigcap_{ \delta >0}  \pi_1(M, \delta))$. Therefore $\tilde{M} / (\bigcap_{ \delta >0}  \pi_1(M, \delta))$ is the universal delta-cover $\tilde{M}^0$.
\qed

\begin{rmrk}\label{univ-delta-cover-pants} 
In Example~\ref{pants}, the universal slipping group of
the surface, $M^2$, was the 
whole fundamental group.  So the universal delta cover 
of $M^2$ is just $M^2$.   Thus all delta covers of this space
are just $M^2$.   
\end{rmrk}

An immediate consequence of Theorem~\ref{univcov}
combined with Proposition~\ref{univ-del-cov-inf} is:

\begin{coro} \label{coro-univ-cov}
If the slipping group is empty
and the covering spectrum has a positive infimum then
the universal cover is a $\delta$ cover.
\end{coro}

\begin{rmrk}\label{open-univ-cov}
If one extends the notion of slipping group and universal
delta cover to complete length spaces which are not known to
have a universal cover, a theorem similar to 
Corollary~\ref{coro-univ-cov} might be applied to prove the existence
of a universal cover.
\end{rmrk}

\section{The Rescaled Covering Spectra}

In this section, we define two scale invariant
spectra called the {\em Rescaled Covering Spectrum}
and the {\em Rescaled Covering Spectrum at Infinity}.
The first of these will be defined for pointed spaces
$(X,x)$ and the latter will not depend on a basepoint.
Both will be invariant when the space is rescaled.

Unlike the ordinary covering spectrum, the rescaled
spectra will be initially defined without covering spaces.
Instead we will assume the spaces have universal covers and
use subgroups of the fundamental group to define the spectra
just as the covering spectrum may be computed using
(\ref{covspecdef2}).
At the end of this section we propose a possible means of extending
the definition to spaces without universal covers
[Remark~\ref{rmrk-extend-rescaled-cov-spec}].  However
we believe the main applications of these spectra are to
submanifolds of Euclidean space and to Riemannian
manifolds both of which always have universal
covers.

\subsection{The rescaled length and the infinite rescaled length}

We begin with the rescaled length:

\begin{defn} \label{defn-rescaled-length}
Given a pointed space $(X,x_0)$ with a universal cover,
we defined the rescaled length of an element $g\in \pi_1(X)$
envisioning it as a deck transform:
\be \label{rescale0}
L^{x_0}_{rs}(g)=
\inf_{\,x \in X\setminus \{x_0\}\,}
\frac{d_{\tilde{X}}(g\tilde{x},\tilde{x})}{\,d_X(x,x_0)\,}.
\ee
where $\tilde{x}$ is a lift of $x$ to the universal cover.
\end{defn}

By definition this is scale invariant, so if we rescale the
metric on our space, we get the same rescaled lengths
for all $g$.

\begin{rmk}\label{supinfrs}
As in Remark~\ref{supinf}, one might
try to extend this definition to spaces without universal covers.
Here, however, we will keep things simple.
\end{rmk}

Note that for a standard cylinder the rescaled length
of all elements are zero, while for the one sheeted hyperboloid,
$\{(x,y,z):\,x^2+y^2=z^2+1\}$, the rescaled length
of the generator is $2$.  See
Example~\ref{one-sheet} for more details.
In both of these examples the infimum is not achieved.

However, for some spaces, the rescaled length is achieved and is
highly dependant on the basepoint.  This occurs with spaces like
the catenoid and those with handles where for most points in
the manifold the shortest geodesic based at a point must
traverse all the way into some central location
and then go all the way back out:

\begin{example} \label{ex-handlebody-rescaled}
Let $X$ be a handlebody with a cusp and $g$ an element with
length, $L(g)<2$, represented by a loop running around the handle.
Suppose $\sigma$ is the minimal loop running around the handle
so $L(\sigma)=L(g)$.  For $x$ far
from the handle, a loop freely homotopic to $\sigma$ based
at $x$ would have to traverse all the way to the loop and then come
back so
\be
\frac{ d_{\tilde{X}} (\tilde{x}, g\tilde{x}) }{ d_X(x,\sigma) } \to 2.
\ee
No matter where $x_0$ is located
\be
\lim_{x\to\infty} \frac{ d_{\tilde{X}}(\tilde{x},g\tilde{x}) }{ d_X(x,x_0) } \to 2.
\ee
On the other hand if we choose $x$ on $\sigma$, we see that
\be
L^{x_0}_{rs}(g)\le
\frac{d_{\tilde{X}}(g\tilde{\sigma(t)},\tilde{\sigma(t)})}
                  {\,d_X(\sigma(t),x_0)\,}
\le
\frac{L(\sigma)}{\, \inf_{t\in {\Bbb{S}}^1}d_X(\sigma(t),x_0)\,}.
\ee
So if we choose a basepoint $x_0$ sufficiently far from
$\sigma$ we get
\be
L^{x_0}_{rs}(g)<2
\ee
So $L^{x_0}_{rs}(g)$ is achieved at some point $x$
rather than approached as $x\to\infty$.  If we
fix a particular such $x_0$ we get some positive value $L^{x_0}_{rs}(g)=f(x_0)>0$.

However this function $f$ is highly dependant on the basepoint
$x_0$ and in fact decays to 0 as we choose $x_0$ further and
further from $\sigma$:
\be
f(x_0)\le \frac{d_{\tilde{X}}(g\tilde{\sigma(t)},\tilde{\sigma(t)})}
                     {\,d_X(\sigma(t),x_0)\,}
\le \frac{L(\sigma)}{\,d_X(\sigma(t),x_0)\,}.
\ee
\end{example}

In our next definition, we define a scale invariant length
which we will prove does not depend on a basepoint
in Lemma~\ref{nodepend} below:

\begin{defn} \label{defn-infinite-rescaled-length}
Given a complete length space $X$ with a universal cover
the infinite rescaled length of an element of the deck transforms is:
\be \label{rescale1}
L_{rs}^\infty(g)=
\lim_{R\to\infty}\inf_{\,x\in X\setminus B_{x_0}(R)\,}
\frac{d_{\tilde{X}}(g\tilde{x},\tilde{x})}{\,d_X(x,x_0)\,}.
\ee
where $\tilde{x}$ is a lift of $x$ to the universal cover
and $x_0 \in X$.
\end{defn}

Here the $rs$ in $L^{x_0}_{rs}$ refer to the word ``rescaled".
They are not parameters.

\begin{lemma} \label{nodepend}
The infinite rescaled length of a deck transform does not depend on
the basepoint $x_0$ used in (\ref{rescale1}).
\end{lemma}

\Pf
Note that there exists $x_i \to \infty$ such that
\be
L_{rs}^\infty(g, x_0)=
lim_{i\to\infty}
\frac{d_{\tilde{M}}(g\tilde{x}_i,\tilde{x}_i)}{d(x_i,x_0)}.
\ee
Let $d(x_0,y_0)=r$.  Then by the triangle inequality,
\be
\frac{d_{\tilde{M}}(x_i, gx_i)}{d(x,y_0)}
\ge \frac{d_{\tilde{M}}(x_i, gx_i)}{d(x,x_0)+ r}.
\ee
Taking $x_i$ to infinity we get:
\be
L^\infty_{rs}(g,y_0) \le
lim_{i\to\infty}
\frac{d_{\tilde{M}}(g\tilde{x}_i,\tilde{x}_i)}{d(x_i,y_0)}
\le L^\infty_{rs}(g,x_0).
\ee
\qed

\begin{lemma}\label{greaterscaling}
Given $X$ and choosing any basepoint $x_0$ to define $L^{x_0}_{rs}$ we have
\be
L^\infty_{rs}(g) \ge L^{x_0}_{rs}(g).
\ee
\end{lemma}

\Pf
This is just a matter of noting that the
the infimum in the definition of $L^\infty_{rs}$ is over a smaller set
than the one in the definition $L^{x_0}_{rs}$.
\qed

\begin{lemma} \label{lessthan2}
The rescaled length and infinite rescaled lengths
are always $\le 2$.
\end{lemma}

\Pf
Let $(X,x_0)$ be a complete length space with
a universal cover and $g \in \pi_1(X)$.
Let $\gamma$ run
between $\tilde{x}_0$ and $g\tilde{x}_0$
in $\tilde{X}$.  Then for any sequence
$x_i \to \infty$ we have
\begin{eqnarray}
L^\infty_{rs}(g) & \le &
lim_{i\to\infty}
\frac{d_{\tilde{M}}(g\tilde{x}_i,\tilde{x}_i)}{d(x_i,x_0)} \\
&\le & lim_{i\to\infty}
\frac{2 d_{\tilde{M}}(\tilde{x}_i,\tilde{x}_0)+ L(\gamma)}{d(x_i,x_0)}
\,=\, 2.
\end{eqnarray}
\qed

\begin{lemma}\label{slipsame}
For a given element $g\in \pi_1(X)$,
$L^{x_0}_{rs}(g)=0$ iff $L_{rs}^\infty(g)=0$.
\end{lemma}

\Pf
If $L^\infty_{rs}(g)=0$, then there exists $x_i \to \infty$
such that
\be
\lim_{i\to\infty}d(\tilde{x_i},g\tilde{x_i})/d(x_0,x_i)=0,
\ee
so the infimum in the definition of $L^{x_0}_{rs}(g)$ is $0$ as well.
We get the converse using (\ref{greaterscaling}).
\qed

\subsection{Rescaled groups and rescaled delta covers}

\begin{defn}\label{defn-rescaled-slip-group}
The rescaled slipping group is the group generated by
$g\in \pi_1(X)$ with $L^\infty_{rs}(g)=0$ and is denoted
\be
\pi^\infty_{rs}(X,0).
\ee
\end{defn}

\begin{defn} \label{defn-rescaled-delta-group}
The rescaled $\delta$ covering group, $\pi^{x_0}_{rs}(X,\delta)$,
of a pointed complete length space $(X,x_0)$ with a universal
cover is the subgroup of $\pi_1(X)$ generated by
elements of rescaled length $L^{x_0}_{rs}(g)<2\delta$:
\be
\pi^{x_0}_{rs}(X,\delta)\,\,=\,\,\lp \, g: \,
L^{x_0}_{rs}(g)<2\delta\, \rp \,\, \subset \,\, \pi_1(X)
\ee
Similarly the infinitely rescaled $\delta$ covering group is
\be
\pi^\infty_{rs}(X,\delta)\,\,=\,\,
\lp g: \,L^\infty_{rs}(g)<2\delta\, \rp  \,\, \subset\,\, \pi_1(X)
\ee
\end{defn}

These groups are scale invariant:
\be \label{groupsscale}
\pi^{x_0}_{rs}(X,\delta)\,=\pi^{x_0}_{rs}(X/R,\delta)\,
\textrm{ and }
\pi^\infty_{rs}(X,\delta)\,=\pi^\infty_{rs}(X/R,\delta).
\ee

\begin{defn} \label{defn-rescaled-delta-cover}
A rescaled $\delta$ cover of a space $X$ with a universal
cover is
\be
\tilde{X}^{\delta,x_0}_{rs} / \pi^{x_0}_{rs}(X,\delta)
\ee
and the infinite rescaled $delta$ cover is
\be
\tilde{X}^{\delta,\infty}_{rs} / \pi^\infty_{rs}(X,\delta)
\ee
\end{defn}

\begin{lemma}\label{lem-rescaled-delta-cover}
We have
\be
\pi^\infty_{rs}(X,\delta) \subset \pi^{x_0}_{rs}(X,\delta).
\ee
Thus
\be
\tilde{X}^{\delta,x_0}_{rs} \textrm{ covers } \tilde{X}^{\delta,\infty}_{rs}.
\ee
When $\delta>1$, the group is the entire fundamental group
and the covering spaces are just the original space $X$.
\end{lemma}

\Pf
Given $g\in \pi^\infty_{rs}(X,\delta)$, $L^\infty_{rs}(g)\le 2\delta$,
so by (\ref{greaterscaling}), $L^{x_0}_{rs}(g)\le 2\delta$,
and $g\in \pi^{x_0}_{rs}(X,\delta)$.
Lemma~\ref{lessthan2} justifies the claims regarding $\delta=1$.
\qed

\subsection{Rescaled covering spectra defined}

\begin{defn}\label{rescaledcovspec}
 Given a pointed complete length space $(X, x_0)$, with a universal
cover, $\tilde{X}$, the rescaled covering spectrum of $X$
denoted $CovSpec_{rs}^{x_0}(X)$,
is the set of all $\delta > 0$ such
 that
 \be
 \pi^{x_0}_{rs}(X,\delta) \neq \pi^{x_0}_{rs}(X,\delta') \qquad \forall \delta'>\delta
 \ee
when viewed as subsets of $\pi_1(X)$.
The infinite rescaled covering spectrum, $CovSpec_{rs}^\infty(X)$,
is defined similarly as the set of all $\delta > 0$ such
 that
 \be
 \pi^\infty_{rs}(X,\delta) \neq \pi^\infty_{rs}(X,\delta') \qquad \forall \delta'>\delta
 \ee
so that it does not depend on the basepoint $x_0$.
\end{defn}

By (\ref{groupsscale}), we see that these covering spectra
are scale invariant:
\be \label{covscale}
CovSpec_{rs}^{x_0}(X)=CovSpec_{rs}^{x_0}(X/R) \textrm{ and }
CovSpec^\infty_{rs}(X)=CovSpec^\infty_{rs}(X/R)
\ee

\begin{prop} \label{covspec-under-1}
 Given a pointed complete length space $(X, x_0)$, with a universal
cover, $\tilde{X}$, 
\begin{eqnarray}
CovSpec_{rs}^{x_0}(X)&\subset& (0,1] \\
CovSpec_{rs}^\infty (X)&\subset& (0,1].
\end{eqnarray}
\end{prop}

\Pf
If $\delta>0 \in CovSpec_{rs}^{x_0} (X)$ then there exists
$g_j \in \pi_1(X) \setminus  \pi^{x_0}_{rs}(X,\delta)$ such
that 
\be
\limsup_{j\to\infty} L^{x_0}_{rs}(g_j)/2 =\delta.
\ee   
By Lemma~\ref{lessthan2} we know $L^\infty_{rs}(g_j)\le 2$.
The proof for infinite rescaled cover spectrum is the same.
\qed

\begin{rmrk}\label{rmrk-semiclosed}
By the theorem in the appendix of \cite{SoWei4} it is easy to
see the rescaled and infinite rescaled covering spectra are 
lower semiclosed sets.
\end{rmrk}

\subsection{Further directions}

\begin{rmrk} \label{rmrk-extend-rescaled-cov-spec}
One should be able to extend the notion of the rescaled covering
spectra to complete noncompact length spaces without universal
covers by adapting the definition of the rescaled lengths
as discussed in Remark~\ref{supinfrs}.
\end{rmrk}

\begin{rmk} \label{rsopen1}
An alternative way of extending the definition of the rescaled
covering spectra to spaces without universal covers
would be to describe the related covering spaces as Spanier covers
with some well chosen sets.  One cannot just choose balls
$B_x(r_x)$ where $r_x=\delta d(x,x_0)$. as there is a difficulty
with the point $x_0$ itself.  If one uses an arbitrary ball around
$x_0$ we lose scale invariance.  One possibility might be to take
a limit of covers where the radius of the ball
about $x_0$ is taken to $0$.
\end{rmk}

\begin{rmk} \label{rsopen2}
To define the infinite rescaled covering spectra on spaces
without universal covers, one needs to figure out how to use
open sets to force the length out to infinity.
\end{rmk}

\section{Asymptotic Behavior and the Rescaled Covering Spectrum}

Here we study the asymptotic behavior of complete noncompact
metric spaces using the rescaled covering spectra and the
rescaled slipping group.  The first subsection concerns the loops to
infinity property and the relationship between the rescaled covering
spectra and the cut off covering spectrum defined by the authors in
\cite{SoWei4}.   The second subsection explicitly computes
the rescaled covering spectra of spectific spaces including hyperboloids and
doubly warped products.   The final subsection proposes a conjecture
concerning the relationship between the diameter growth of a space
and its rescaled covering spectrum.

\subsection{Loops to infinity}

We recall the loops to infinity property defined in
\cite{So-loops}:

\begin{defn}\label{defloopstoinfty}
Given a metric space, $X$,
a loop $\gamma:{\Bbb{S}}^1 \to X$ is said to have the {\em
loops to infinity} property, if for every compact set $K \subset X$,
there is another loop $\sigma:{\Bbb{S}}^1 \to X\setminus K$
freely homotopic to $\gamma$.
\end{defn}

\begin{lemma} \label{rsloopstoinfty}.
If $L_{rs}^\infty(g)<2$ then any curve
representing $g$ has the loops to infinity property.
\end{lemma}

\Pf
If $g$ does not have the loops to infinity property,
its representative curves $C_i$ whose lengths approach $L^{\infty}_{rs}(g)$
and have $C_i(0)\to\infty$ must pass through a common compact
set, $K$.  Thus $L(C_i) \ge 2 d_X(C_i(0),K)$, so the
limit in the definition of rescaled length gives a $2$.
\qed

\begin{theo} \label{rescalevscut}
Given a complete length space $X$, with a universal
cover, $\tilde{X}$, if
\be
CovSpec^\infty_{rs}(X) \in (0,1)
\ee
then
\be
CovSpec_{cut}(X)=\emptyset.
\ee
\end{theo}

This theorem captures the fact that the rescaled covering spectrum and
the cut off spectrum measure very different kinds of ``holes"
in a complete length space.

\Pf
By the hypothesis we have
\be
 \pi^\infty_{rs}(X,1) =\pi^\infty_{rs}(X,\delta') \qquad \forall \delta'>1.
\ee
So for any $g_0\in \pi_1(M)$, taking $\delta'>L(g_0)/2$, we see that
$g_0 \in \pi^\infty_{rs}(X,1)$.   So
the fundamental group is generated
by elements $g$ with $L^{x_0}_{rs}(g)<2$.  Thus by
Lemma~\ref{rsloopstoinfty}, any curve representing such a generator has
the loops to infinity property and so by Theorem 4.20 in \cite{SoWei4} $CovSpec_{cut}(X)=\emptyset
$.
\qed

\subsection{Cones, hyperboloids and warped products}

\begin{defn}
The base point free cone over $Y$ with scaling $k$
denoted $C_k(Y)$ is 
\be
C_k(Y) = (0,\infty)_f \times Y  \label{C_k(Y)}
\ee
where $f(r)=kr$.
\end{defn}


\begin{theo} \label{thmcone}
Given a compact length space $Y$, we have
\be
CovSpec^\infty_{rs}(C_k(Y))=
\Big\{ (1/2)\sqrt{ \, 2 -2 \cos(\,\min\{\pi, 2k \delta\} ) }\,\,: \,\, \delta \in CovSpec(Y)\Big\}. 
\ee
and
\be
CovSpec^{x_0}_{rs}(C_k(Y))=\emptyset.
\ee
\end{theo}

\Pf
Recall that with a linear
warped product, 
\be
d_{C_k(Y)}((y_1, r_1), (y_2, r_2))=
\sqrt{ \, r_1^2 + r_2^2 -2r_1 r_2 \cos(\, \min\{\pi,kd_Y(y_1, y_2) \})\,}.
\ee  
This can be seen because if $kd_Y(y_1, y_2)<\pi$ then
the minimal geodesic between them has length $<\pi$,
and the minimal geodesic between $(y_i,r_i)$ lies
in a linear warped product of that geodesic (which is a sector
of Euclidean space and one can compute its length using the
law of cosines.  If $kd_Y(y_1,y_2)\ge  \pi$ then the
shortest geodesic in the cone passes though the base point
and has length 
\be
r_1+r_2=\sqrt{ \, r_1^2 + r_2^2 -2r_1 r_2 \cos(\pi)\,}.
\ee  

Let $N=C_k(Y)$ and $\tilde{N}$ its universal cover.
By the definition of the base point free cone,
any $g\in \pi_1(N)$ is also a $g\in \pi_1(Y)$.
In fact, one can easily see that 
\be
\tilde{N}= (0,\infty) \times_f \tilde{Y}=C_k(\tilde{Y}),
\ee
where any point in $N$ can be represented as 
$x=(y,r)$ with $x\in Y$ and this point lifts to
a point $\tilde{x}=(\tilde{y}, r) \in \tilde{N}$ where
$\tilde{y}\in \tilde{Y}$ and where $g(\tilde{x},r)=(g\tilde{x},r)$.

Let $y_1$ achieve the infimum in the definition of $L(g)$
viewing $g\in \pi_1(Y)$:
\be
L(g)=d_{\tilde{Y}}(g\tilde{y_1}, \tilde{y_1})
\ee
and let $y_0$ be furthest from $y_1$:
\be
d_Y(y_0, y_1)=\sup_{y\in Y} d_Y(y_1,y).
\ee

For any $r_0>0$ we take $x_0=(y_0, r_0)$.   We have
\be
r^{-1}(0,R) \subset B_{x_0}(R+r_0) \subset r^{-1}(0, R+2r_0).
\ee
Applying  Definition~\ref{defn-rescaled-length} which
does not depend on the choice of $x_0$, we have 
\begin{eqnarray} 
L_{rs}^\infty(g) &=& \lim_{R\to\infty} 
\inf_{\,(y,r) \in N\setminus B_{x_0}(R)\,}
\frac{d_{\tilde{N}}(g(\tilde{y},r),(\tilde{y},r))\,}{\,d_N((y,r),(y_0,r_0))\,} \\
&=&\lim_{R\to\infty} 
\inf_{y \in N, r\ge R\,}
\frac{d_{\tilde{N}}((g\tilde{y},r),(\tilde{y},r))\,}{\,d_N((y,r),(y_0,r_0))\,} \\
&=&\liminf_{r\to\infty} 
\inf_{y \in N\,}
\frac{r \sqrt{ \, 2 -2 \cos(\,\min\{\pi, kd_{\tilde{Y}}(g\tilde{y}, \tilde{y})\} )\,}
\,}{\,
\sqrt{ \, r^2 + r_0^2 -2r r_0 \cos(\, \min\{\pi, kd_Y(y, y_0)\} )\,}
\,} \\
&=&\liminf_{r\to\infty} 
\frac{r \sqrt{ \, 2 -2 \cos(\,\min\{\pi, kd_{\tilde{Y}}(g\tilde{y}_1, \tilde{y}_1)\} )\,}
\,}{\,
\sqrt{ \, r^2 + r_0^2 -2r r_0 \cos(\, \min\{\pi, kd_Y(y_1, y_0)\} )\,}
\,} \\
&=&
 \sqrt{ \, 2 -2 \cos(\,\min\{\pi, kd_{\tilde{Y}}(g\tilde{y}_1, \tilde{y}_1)\} )\,} \\
 &=&
 \sqrt{ \, 2 -2 \cos(\,\min\{\pi, k L(g)\} )\,} .
\end{eqnarray}
This implies the first claim in the statement of our theorem.

In contrast, when we compute the rescaled length depending on a given
basepoint $x_0=(y_0,r_0)$, we choose any $y_2\neq y_0$
and we have
\begin{eqnarray} 
L_{rs}^{x_0}(g) &=& 
\inf_{\,(y,r) \neq (y_0,r_0)\,}
\frac{d_{\tilde{N}}(g(\tilde{y},r),(\tilde{y},r))\,}{\,d_N((y,r),(y_0,r_0))\,} \\
&=&
\inf_{\,(y,r) \neq (y_0,r_0)\,}
\frac{r \sqrt{ \, 2 -2 \cos(\,\min\{\pi, kd_{\tilde{Y}}(g\tilde{y}, \tilde{y})\} )\,}
\,}{\,
\sqrt{ \, r^2 + r_0^2 -2r r_0 \cos(\, \min\{\pi, kd_Y(y, y_0)\} )\,}
\,} \\
&\le&
\inf_{\,r\neq r_0\,}
\frac{r \sqrt{ \, 2 -2 \cos(\,\min\{\pi, kd_{\tilde{Y}}(g\tilde{y}_2, \tilde{y}_2)\} )\,}
\,}{\,
\sqrt{ \, r^2 + r_0^2 -2r r_0 \cos(\, \min\{\pi, kd_Y(y_2, y_0)\} )\,}
\,} =0.
\end{eqnarray}
\qed

This effect where the rescaled spectrum goes to $0$ does not
appear to happen in a setting where the space $N$ is a
manifold that is asymptotically 
cone like. Although the infinite rescaled spectrum 
appears to behave in the
same way.   This can be seen in the following example: 

\begin{example}\label{one-sheet}
The one-sheeted hyperboloid, 
\be
N^2=\{(x,y,z):x^2+y^2=z^2+1\},
\ee
has a rescaled covering spectrum which consists
of a single value
\begin{eqnarray}
CovSpec_{rs}^{x_0}(N)&=&\left\{ 1 \right\}\\
CovSpec_{rs}^\infty(N)&=&\left\{ 1  \right\}.
\end{eqnarray}
\end{example}

\Pf
To compute $L^{x_0}_{rs}(g)$, we
take advantage of the invariance under rescaling:
\begin{eqnarray*}
\{(x,y,z):x^2+y^2=z^2+1\}/R&=&
\{(Rx,Ry,Rz):x^2+y^2=z^2+1\}\\
&=&\{(x,y,z):x^2+y^2=z^2+1/R^2\}
\end{eqnarray*}
Let $\sigma\subset N$ be the circular neck.
As $R\to\infty$, we see that the length of a shortest geodesic $\gamma$
based at any point $\gamma(0)\in N\setminus T_r(\sigma)$
is approaching
the length of such a geodesic in the flat cone with the origin removed:
\be
\{(x,y,z):x^2+y^2=z^2, \, (x,y,z)\neq(0,0,0)\}.
\ee
This is $C_k (Y)$ in (\ref{C_k(Y)})  with $k=1, Y$  the circle with radius $1/\sqrt{2}$.  Covspec$(Y) =\{ \frac{\pi}{\sqrt{2}} \}$, hence Theorem~\ref{thmcone} gives 
\[
CovSpec_{rs}^\infty(C_k(Y)) = \{1\}.\]
Since it is scale invariant, so CovSpec$_{rs}^{\infty}(N^2) = \{1\}$.
\qed

In the following theorem we see that warped product spaces
which are asymptotic to such cones, and even double warped producs
which need not be close to the cones but are simply connected
on their second warping factor, have the same rescaled covering
spectra as their cones.

\begin{theo}\label{asymthm}  
Suppose our space 
is  a doubly warped product
\be \label{asymthm1}
X = [0,\infty) \times_f N \times_h M
\ee
where $N$ and $M$ are compact Riemannian
manifolds
with $f'(0)=0$, $f(0)>0$
and either $h'(0)=0$ or $h(0)=0$.  In all these cases we 
assume $f(r)>0$ and $h(r)>0$ on $(0,\infty)$.

If $M$ is simply connected or $h(0)=0$ we have
\be
L_{rs}^\infty(g) = \liminf_{r\to\infty} \frac{F( r, L_N(g))}{r}.
\ee
where $F(r,d)$ is the length of a minimal geodesic in
\be
\mathbb{R}\times_{f}\mathbb{R}
\ee
between $(r,0)$ and $(r,d)$ and $L_N(g)$ is the length
of $g$ viewed as an element of $\pi_1(N)$.
\end{theo}

\Pf
Observe that the universal cover
\be
\tilde{X}= [0,\infty) \times_f \tilde{N} \times_h M
\ee
where $x=(r,y,z)$ lifts to $\tilde{x}=(r, \tilde{y}, z)$
and any $g\in \pi_1(X)$ can be viewed as
$g\in \pi_1(N)$ so that
\be
g\tilde{x}=(r, g\tilde{y}, z).
\ee
Also we note that if
\be
X'= [0,\infty) \times_f N 
\ee
then $\tilde{X}'=  [0,\infty) \times_f \tilde{N}$ and that
\be
 d_{X'}((r_1, y_1), (r_2, y_2)) 
\le  d_{X}((r_1, y_1, z_1), (r_2, y_2, z_2))
\le  d_{X'}((r_1, y_1), (r_2, y_2)) + \max_{[r_1, r_2]} h(r) d_M(z_1,z_2).
\ee

So applying  Definition~\ref{defn-rescaled-length} with
$x_0=(0, y_0,z_0)$ 
we have
\begin{eqnarray} 
L_{rs}^\infty(g)  &=& \lim_{R\to\infty}
\inf_{\,(r,y,z) \in X\setminus B_R(0,y_0,z_0)\,}
\frac{d_{\tilde{X}}((r,g\tilde{y},z),(r,\tilde{y},z))}{\,d_X((r,y,z),(0, y_0,z_0))\,} \\
&=& \liminf_{r\to\infty}
\inf_{\,y\in N, z\in M }
\frac{d_{\tilde{X}}((r,g\tilde{y},z),(r,\tilde{y},z))}{\,d_X((r,y,z),(0, y,z))
\pm d_X((0,y_0,z_0),(0,y,z))\,} \\
&=& \liminf_{r\to\infty}
\inf_{\,y\in N, z\in M }
\frac{d_{\tilde{X}}((r,g\tilde{y},z),(r,\tilde{y},z))}{\,d_X((r,y,z),(0, y,z))\,} \\
&=& \liminf_{r\to\infty}
\inf_{\,y\in N}
\frac{d_{\tilde{X}'}((r,g\tilde{y}),(r,\tilde{y}))}{\,d_{X'}((r,y),(0, y))\,} \\
&=& \liminf_{r\to\infty}
\inf_{\,y\in N}
\frac{F(r,d_{\tilde{N}}(g\tilde{y},\tilde{y})) }
{\,r\,}.
\end{eqnarray}

Observe that $F(r,d)$ is increasing as $d$ increases
so
\be
\inf_{\,y\in N}
\frac{F(r,d_{\tilde{N}}(g\tilde{y},\tilde{y})) }{\,r\,} =\frac{F( r, L_N(g))}{r}
\ee
where $L_N(g)$ be the length of $g$ as an element of $\pi_1(N)$.
Thus
\be
L_{rs}^\infty(g) = \liminf_{r\to\infty} \frac{F( r, L_N(g))}{r}.
\ee
\qed

\begin{coro} \label{asymcoro}
Suppose that
\be
\lim_{r\to\infty} \frac{f(r)}{r} =k.
\ee
If $k>0$ we have
\be
CovSpec_{rs}^{\infty}(X)= \{\frac{1}{2} \sqrt{2-2\cos(\min\{\pi, 2kd\}): \,\, d\in CovSpec(N) } 
\ee
and the rescaled slipping group is empty.  If $k=0$
then
\be
CovSpec^\infty_{rs}(X)=\emptyset.
\ee
and $\pi_1(X) \subset \pi_{rs}^\infty(X, 0)$. 
\end{coro}

\Pf
In the special case where 
\be
\lim_{r\to\infty}f(r)/r=0,
\ee
then 
\be
\lim_{r\to\infty} F(r,d)/r =0
\ee
so $CovSpec_{rs}^{\infty}(X)=\emptyset$
and $\pi_1(X)$ is a subset of the
rescaled slipping group.
Alternately if
\be
\lim_{r\to\infty}f(r)/r=k>0,
\ee
then 
\be
\lim_{r\to\infty} F(r,d)/r = \sqrt{2-2\cos(\min\{\pi, 2kd\}) } 
\ee
so we have our claim.
\qed

\subsection{Diameter growth}

For any metric space $X$, $p \in X$, 
\be
\limsup_{r\to\infty} \frac{diam(\partial B_p(r))}{r} =a,
\ee
with $0 \le a \le 2$. 

\begin{conj} \label{conj-diam-growth}
For a metric space $X$ with $a$ above, 
\be
CovSpec_{rs}^\infty(X) \subset (0, \frac a2 ] \cup \{1\}.
\ee  
In particular, when $X$ has sublinear diameter growth, i.e. $a=0$, then  
\be
CovSpec_{rs}^\infty(X) \subset \{1\}.
\ee
\end{conj}

\section{Applications with Curvature Bounds}

In this section we restrict ourselves to the study of complete
Riemannian manifolds, $M$.   In the first subsection we study manifolds with nonnegative
sectional curvature, whose topology can be understood because such
manifolds have a compact soul.   We also present an interesting 
example of Wilking which has nonnegative sectional curvature and
positive Ricci curvature and prove it has a nontrivial rescaled 
covering spectrum.   In the next subsection we study manifolds with nonnegative Ricci curvature.   We first show that the fundamental
group of the classic example of Nabonnand lies in the rescaled 
slipping group [Proposition~\ref{warp-rescale}].   We prove Theorem~\ref{finite}
concerning manifolds with positive Ricci curvature and its
Corollary~\ref{coro-nondec}.   We close with a subsection on
the rescaled covering spectra of
manifolds with nonnegative Ricci curvature proving 
Theorem~\ref{Ricci-pi-rs-thm}. 

\subsection{Nonnegative sectional curvature}

Cheeger-Gromoll \cite{Cheeger-Gromoll} proved that complete manifolds 
with nonnegative
sectional curvature are diffeomorphic to normal bundles over totally geodesic
compact submanifolds called souls.  Sharafutdinov \cite{Sharafutdinov} then proved there
was a distance nonincreasing retraction to the soul: $P:M \to S$.
Perelman \cite{Perelman1, Perelman2} showed that 
$P$ is a Riemannian submersion and extended the distance 
nonincreasing retraction to complete Alexandrov spaces 
with nonnegative curvature.  In \cite{SoWei4} using the distance nonincreasing 
retraction, $P$, the authors proved that the covering spectrum of these 
spaces behave exactly like the covering spectrum of a compact space:

\begin{theo} \label{oldsectthm} \cite{SoWei4}
If $M^n$ is a complete noncompact Alexandrov space with nonnegative
curvature, then 
\be
CovSpec(M^n)=CovSpec(S^k)
\ee
 where $S^k$ is its soul and 
 \be
 CovSpec(M^n) \subset (1/2)Length(M^n)= (1/2)Length(S^k)
 \ee
and it is
determined by the marked length spectrum of $M^n$.
\end{theo}

In the proof of this theorem, we observed that the length of
every element of $\pi_1$ was achieved within the soul.  Combining this
with our Theorem~\ref{thm-line}, we see that either there are only finitely
many elements in the fundamental group or the universal
cover contains a line.   This was already shown in the
compact setting by Cheeger-Gromoll \cite{ChGr1971}.   

\begin{theo} \label{newsectthm} \cite{SoWei4}
If $M^n$ is a complete noncompact Alexandrov space with nonnegative
curvature, then immediately
\be
Shift(M^n)=Length(M^n),
\ee
the slipping group and universal slipping groups are trivial.
\end{theo}

We can also prove the following theorem relating the holonomy
to the covering spectrum:

\begin{theo} \label{sect-holonomy}
If  a complete noncompact manifold, $M^m$,
with nonnegative sectional curvature has a connected holonomy group, then 
the fundamental group lies in the rescaled slipping group 
and the infinite rescaled covering spectrum is empty.  More generally, 
if $g\in \pi_1(M^n)$ has a representative loop whose holonomy
lies in the identity component of the holonomy group
then $\sup\{L(g,p): p\in M^m\}<\infty$ and
$g$ lies in the rescaled slipping group.   
\end{theo}

\Pf
Perelman proved that for all $t>0$, we have
\be
d_M(\exp_{q_1}(tv_1), \exp_{q_2}(tv_2)) \ge d_S(q_1, q_2) 
\ee
as long as $v_i\in TM_{q_i}$ are perpendicular to the soul.   Furthermore,
he proved that
if there is a path from $q_1$ to $q_2$ such that $v_1$ is parallel
to $v_2$ along that path, then there is equality for all $t>0$ \cite{Perelman1}\cite{Perelman2}.

Now suppose $g\in \pi_1(M)$ has a representative loop, $C$, which lies
in the identity component of the holonomy group of $M^n$.   
We apply the distance decreasing retraction, $P$, to obtain a loop
$P\circ C$, in the soul which is freely homotopic to $C$, and thus
also has holonomy in the identity component of the holonomy group
of $M^n$.   Recall the identity component
of the holonomy is generated by contractible loops. So in this case one
can compose $P\circ C$ with a trivial loop $\gamma_0$ so that
the composition, $\gamma:[0,1]\to S\subset M$, 
has a closed parallel normal vector field, $v(s)$, such that $v(0)=v(1)$.   
Thus, taking $C_t(s)=\exp_{\gamma(s)}(tv(s))$ as our representative
of $g$ based at $p_t=C_t(0)$ we have
\be
L(g,p_t) \le L(C_t) = L(\gamma).
\ee
This
\be
L_{rs}^\infty(g) \le \liminf_{t\to\infty} \frac{L(g, p_t)}{t} =0.
\ee
\qed

\begin{rmrk}
Note that the infinite Moebius strip has holonomy group, $\mathbb{Z}_2$,
and the value $1$ lies in its infinite rescaled covering spectrum.
One might think that any complete noncompact manifold with
nonnegative sectional curvature that has a nonempty infinite rescaled
covering spectrum has a line in its universal cover.  However,
in the appendix we will show that Wilking's example has a nontrivial
covering spectrum.
\end{rmrk}

\subsection{Positive Ricci curvature and Theorem~\ref{finite}}

For many years after Cheeger and Gromoll proved the Soul Theorem \cite{Cheeger-Gromoll}
for manifolds with nonnegative sectional curvature, it was an open
question whether a complete noncompact manifold with positive
Ricci curvature must have a finite fundamental group.   After all,
if the manifold has positive Ricci curvature and nonnegative 
sectional curvature, then it has a soul with positive Ricci curvature
and the fundamental group must be finite by Myers' Theorem \cite{Myers}.
Manifolds with positive sectional curvature have trivial fundamental
groups.   

Then in 1980, Nabonnand \cite{Nab} found an
example of a complete noncompact manifold with positive Ricci
curvature whose fundamental group is $\mathbb Z$.
In fact, in \cite{Wei1988}, the second author showed that fundamental
group could be any torsion free nilpotent group and Wilking adapted these examples
 to prove the fundamental group could be any almost
nilpotent group \cite{Wilking2000}.

  Below
we will describe Nabonnand's example in more detail and will see
below in Proposition~\ref{warp-rescale}, that all elements
of the fundamental group  of his example lie in the 
rescaled slipping group.   We will also study the Ber\'ard-Bergery
Examples which demonstrate a variety of covering spectra can
be achieved on manifolds with positive Ricci curvature
[Proposition~\ref{prop-covspec-ric}].   In the Appendix we
will see that Wilking's example has positive Ricci curvature and
a nontrivial rescaled covering spectrum.

First we prove Theorem~\ref{finite}, stated in the introduction,
that for any compact set $K$ in a complete noncompact manifold
with positive Ricci curvature, there are only finitely many distinct elements $g$ of $\pi_1(M)$ such 
that  $L(g)$ are achieved in the compact set $K$.

\Pf
By Theorem~\ref{thm-line}, we know that if there are
infinitely many elements whose length is achieved within a compact
set, then the universal cover contains a line.   However, by the
Cheeger-Gromoll Splitting Theorem \cite{Cheeger-Gromoll}, any manifold
with nonnegative Ricci curvature that contains a line splits isometrically,
and thus cannot have strictly positive Ricci curvature.
\qed

A corollary of Theorem~\ref{finite} is:

\begin{coro} \label{coro-nondec}
Let $M$ be a complete Riemannian manifold with $\Ric > 0$. If there is a compact set $K \subset M$ such that $L(g,p)$ is nondecreasing outside of $K$ as $p \to \infty$,  then $\pi_1(M)$ is finite.
\end{coro}

\begin{rmrk}
Note that the nondecreasing hypothesis in Corollary~\ref{coro-nondec}
implies that the slipping group is trivial and the covering
spectrum has an infimum.   One may ask if one still must
have a finite fundamental group when the nondecreasing hypothesis
is replaced by an assumption that the slipping group is trivial
and the covering spectrum has an infimum.   However we will
produce many examples in Proposition~\ref{prop-covspec-ric}
which have infinite fundamental groups. 
\end{rmrk}

Then Ber\'ard-Bergery \cite{Berard-Bergery} showed that, given 
any compact manifold $M^m$ with
$Ricci \ge 0$, $N^{m+3}=M^m\times \mathbb R^3$ has a complete warping metric
\be  \label{warping}
g = dr^2 + h^2(r) g_{{\Bbb{S}}^2}  + f^2(r) g_M \textrm{ with } h(0)=0,
\,\, h'(0)=1,
\textrm{ and } f(0)\not= 0, f'(0) =0.
\ee
with $\Ric >0$.  Nabonnand's example had the same structure
with $M={\Bbb{S}}^1$.  

For complete manifolds with $\Ric \ge 0$, unlike manifolds with 
nonnegative sectional curvature, $L(g,p)$ could  decrease as $p$ goes 
to infinity. In Nabannand's example \cite{Nab} of 
${\Bbb{S}}^1 \times \mathbb R^3$ with $Ric >0$, the length of ${\Bbb{S}}^1$ 
strictly decreases as it goes to infinity.  In fact we can show that:

\begin{prop} \label{warp-rescale}
For any warping metric as in (\ref{warping}) (with $M = {\Bbb{S}}^1$) on 
$N={\Bbb{S}}^1 \times \Bbb{R}^3$), if $Ric >0$, then $f'(r) < 0$ for $r>0$.
So the rescaled covering spectrum is empty
and $\pi_1(N)$ is contained in the rescaled slipping group.
\end{prop}

\Pf
To see this we examine the Ricci curvature in the ${\Bbb{S}}^1$ direction, $V$,
which, according to Nabonnand's equation (1) is
\be \label{nabeq}
- Ric(V,V)=\frac{f''(r)}{f(r)}  + 2 \frac{f'(r)h'(r)}{f(r)h(r)} < 0
\qquad{ for }\ r>0.
\ee
By smoothness requirements $f'(0)=0$ and $f(r)>0$ for all $r\in [0,\infty)$.  
Suppose $f'(r_0)=0$.  Then
\be
\frac{f''(r_0)}{f(r_0)}  =\frac{f''(r_0)}{f(r_0)}  
+ 2 \frac{f'(r_0)h'(r_0)}{f(r_0)h(r_0)} < 0.
\ee
So $f''(r_0) <0$ and all critical points are local maximal.  This implies there is at most
one critical point and $f'(r)<0$ after that critical point.
Since $f'(0)=0$ we have $f'(r)<0$ for all $r >0$. So
$f$ is strictly decreasing.

Any $g\in \pi_1(N)$
is also a $g\in \pi_1(M)$ with such a construction
and any point in $N$ can be represented as $(x,r)$ with $x\in M$.  
So applying  Definition~\ref{defn-rescaled-length}
\begin{eqnarray} 
L_{rs}^\infty(g) \le L^{x_0}_{rs}(g) &=&
\inf_{\,(x,r) \in N\setminus \{(x_0,0)\}\,}
\frac{d_{\tilde{N}}(g\tilde{(x,r)},\tilde{(x,r)})}{\,d_N((x,r),(x_0,0))\,} \\
&\le&
\lim_{r\to\infty}\frac{f(r)d_{\tilde{M}}(g\tilde{x},\tilde{x})}
{\, r -d_N(x,x_0)\,}
=0.
\end{eqnarray}
\qed

\begin{rmrk}
The above theorem also holds for general warped products as in 
(\ref{warping}) when
 $M^n$ has nonnegative Ricci curvature, e.g. $M^n = T^n$. 
 \end{rmrk}
 
 \begin{rmk}  
Note that in general the rescaled slipping group of a manifold
with $Ricci>0$ may not be the entire fundamental group
even though such manifolds have the
loops to infinity property.  See Example~\ref{wilking} in the
appendix.
\end{rmk}
 
\begin{rmrk} 
In Ber\'ard-Bergery type examples
where $f(r)$ is monotone decreasing in $r$ in (\ref{warping}),
$\lim_{r \ra \infty} f(r)$ could be zero or any positive number. 
When $\lim_{r \ra \infty} f(r) = 0$, the entire fundamental
group is in the slipping group and the universal delta cover is the 
whole space.    
\end{rmrk}

\begin{prop}\label{prop-covspec-ric}
If $M^n$ is a compact manifold with nonnegative
Ricci curvature, then there is a Ber\'ard-Bergery example, $N^{n+3}$
of a complete noncompact manifold with positive Ricci curvature
as in (\ref{warping}) such that
\be
CovSpec(N^{n+3})= CovSpec(M^n).
\ee
\end{prop}

\Pf
In the Ber\'ard-Bergery examples, $f(r)$ need only be decreasing and
we can find decreasing $f(r)$ such that  $\lim_{r \ra \infty} f(r)=1$.
\qed

 \begin{rmrk}
  If $M^n$ has positive Ricci curvature, then $f$ need not be strictly decreasing to make $M\times \mathbb R^3$ have positive Ricci curvature. One can take $f(r)$ to be constant and $h(r)$ concave down, since the isometric
product of  $M$ with $\Bbb{R}^3$ where $\Bbb{R}^3$ has 
$Ric > 0$ would have positive Ricci curvature.  Then we 
get 
$CovSpec(M\times \Bbb{R}^3)=CovSpec(M)$ and 
$CovSpec_{rs}^{x_0}(M\times \Bbb{R}^3)=0$.  
\end{rmrk}

\subsection{Ricci curvature and the rescaled covering spectrum}

We now prove Theorem~\ref{Ricci-pi-rs-thm-2}.  In fact we prove
the following stronger theorem:

\begin{theo}\label{Ricci-pi-rs-thm-2}
Suppose $M^n$ is a complete noncompact manifold, $M^n$, 
and $\pi_{rs}^{\infty}(M,\delta)$ is trivial.  
If $M^n$ has nonnegative Ricci curvature 
then the fundamental group is finite with less than $2(2+\delta)^n/\delta^n$ elements.    
\end{theo}

\Pf
Assume on the contrary there are $N$ nontrivial
elements $g_1, g_2, ... g_{N-1}\in \pi_1(M)$ 
such that $g_i \neq g_j^{\pm 1}$ for all $i\neq j$.
Let $g_0$ be the identity in $\pi_1(M)$.
When $M$ has nonnegative Ricci curvature, we have our
theorem
if we prove $N \le \frac{(2+\delta)^n}{\delta^n}$.
So we assume on the contrary that  
\be
N > \frac{(2+\delta)^n}{\delta^n}.
\ee

For all $i,j \in \{0,1,2,...(N-1)\}$ with $i \not= j$, we have $g_ig_j^{-1}$ are nontrivial.
So 
\be
g_ig_j^{-1} \notin \pi_{rs}^\infty(M, \delta)
\ee
and
\be
L^\infty_{rs}(g_i g_j^{ -1}) \ge 2\delta.
\ee
So for any $x_0\in M$, $\epsilon >0$, there exists $R_{i,j,\vare}$ sufficiently large that 
$\forall R \ge R_{i,j,\vare}$ we have
\be \label{inf-dist}
\inf_{x\in M\setminus B_{x_0}(R)} 
\frac{d_{\tilde{M}}(\tilde{x}, g_ig_j^{-1}\tilde{x} )}{d_M(x,x_0)} > 2\delta-2\vare.
\ee
Let $R_\vare=\max\{R_{i,j,\vare}: \,\, i,j \in 0,1,...(N-1)\}$.   
Then (\ref{inf-dist}) holds for all $R\ge R_\vare$. 

In particular, for all $x\in M\setminus B_{x_0}(R_\vare)$ we
have
\be
B_{g_i\tilde{x}}(\rho_\vare(x)) \cap B_{g_j\tilde{x}}(\rho_\vare(x)) =\emptyset
\ee
where 
\be
\rho_\vare(x)= \rho(x)(\delta-\vare) \textrm{ where }\rho(x)=d(x,x_0). 
\ee
Fix a lift of $x_0$, $\tilde{x}_0$,  choose $\tilde{x}$,  a lift of $x$,  such that
\be
d_{\tilde{M}}(\tilde{x}, \tilde{x}_0)=d_M(x, x_0)=\rho(x).
\ee

We now apply the trick in Milnor's paper \cite{Mi}.  First
we set
\be
C= \max\{ d_{\tilde{M}}(\tilde{x}_0, g_i \tilde{x}_0):\,\, i=0,1,2...N \}.
\ee
Then by the triangle inequality,
\be
B_{g_i\tilde{x}}(\rho_\vare(x)) 
\subset B_{\tilde{x}_0}(C + \rho(x)+\rho_\vare(x))
\subset B_{\tilde{x}}(C + 2\rho(x)+\rho_\vare(x)).
\ee
Since
\be
\vol(B_{g_i\tilde{x}}(\rho_\vare(x)))=\vol(B_{\tilde{x}}(\rho_\vare(x))  
\ee
and the balls are disjoint, we have
\be
N\vol(B_{\tilde{x}}(\rho_\vare(x))   \le 
\vol(B_{\tilde{x}}(C + 2\rho(x)+\rho_\vare(x)).
\ee
By the Bishop-Gromov Volume Comparison Theorem, we have
\be
N\le
\frac{(C + 2\rho(x)+\rho_\vare(x))^n}
{(\rho_\vare(x))^n}.
\ee
Substituting $\rho_\vare(x)=\rho(x)(\delta-\vare)$ we have
\be
N\le
\frac{(C + \rho(x)(2+\delta-\vare))^n}
{(\rho(x)(\delta-\vare))^n}
\ee
for all $x\in M\setminus B_{x_0}(R_\vare)$.
So now we take $x\to \infty$, $\rho(x)\to \infty$ and we have
\be
N\le
\frac{(2+\delta-\vare)^n}
{(\delta-\vare)^n}.
\ee
Lastly we take $\vare \to 0$ and we have our contradiction.
\qed

\begin{rmrk} \label{lim-Ricci-pi-rs-thm}
A version of Theorem~\ref{Ricci-pi-rs-thm} holds for any metric measure space whose
universal cover satisfies a doubling condition for arbitrarily
large radii including the pointed Gromov-Hausdorff
limits of manifolds with nonnegative Ricci curvature.   In that
case we can prove the group of deck transforms of the universal
covering space is finite.  See \cite{SoWei1} for similar discussions
of extensions of theorems regarding the universal covers of
manifolds with nonnegative Ricci curvature to their limit spaces.
\end{rmrk}

\begin{rmrk}
It should be noted that Anderson
extended Milnor's techniques by intersecting balls with fundamental
domains and taking their radii to infinity in \cite{Anderson}.
He proved that the
fundamental group in a manifold with nonnegative
Ricci curvature and $\liminf_{R\to\infty}\vol(B_{x_0}(R))/R^n =\alpha>0$
has a fundamental group which is finite with
$\le \omega_n/\alpha$ elements.   Li also proved this
using the heat kernal in \cite{Li1986}.  One might ask if
manifolds satisfying the conditions in their theorem have
a trivial rescaled covering group for some value $\delta$.
Or perhaps, one could show that at least the rescaled
slipping group is trivial.
\end{rmrk}

\begin{rmrk}
The first author proved in \cite{Sormani-group} that the fundamental group
of a manifold, $M^n$, with nonnegative Ricci curvature
and linear volume growth
\be
\limsup_{R\to\infty} \frac{\vol(B_{x_0}(R))}{R} <\infty
\ee
has a finitely generated fundamental group.   
One might ask whether under these conditions
one can prove there is a $\delta>0$ such that
$\pi_{rs}^\infty(M, \delta)$ is trivial and thus improve this
theorem to imply that the fundamental group is finite.
However the cylinder and the infinite Moebius strip
both have linear volume growth (and 0 curvature),
and infinite fundamental group with a nontrivial
rescaled slipping group and so $\pi_{rs}^\infty(\delta)$
is nontrivial for all $\delta>0$.
\end{rmrk}

\begin{rmrk}
The first author proved in \cite{Sormani-group} that the fundamental group
of a manifold, $M^n$, with nonnegative Ricci curvature
with small linear diameter growth 
\be
\limsup_{R\to\infty} \frac{\diam(B_{x_0}(R))}{R} \le S_n
\ee
has a finitely generated fundamental group.   The constant
provided in \cite{Sormani-group} was improved by
\cite{XWY}.   One might ask whether under these conditions
one can prove there is a $\delta>0$ such that
$\pi_{rs}^\infty(M, \delta)$ is trivial and thus improve this
theorem to imply that the fundamental group is finite.
However,  this is not possible: Nabonnand's example has both an infinite
fundamental group and sublinear diameter growth
as well as a nontrivial rescaled slipping group.   
\end{rmrk}

\sect{Appendix}

We would like to thank Burkhard Wilking for suggesting this example
when asked for an example of a manifold with nonnegative sectional
curvature that had $L(g,p)$ diverging to infinity as $p\to \infty$,
but the universal cover has no line.

\begin{ex}  \label{wilking}
Let $M^6 = ( {\Bbb{S}}^3 \times \Bbb{R}^4 )/ \textrm{Pin}(2)$. This has nonnegative
sectional curvature because it is a quotient of a manifold with 
nonnegative sectional curvature. It is
$\mathbb R^4$ bundle with soul ${\Bbb{S}}^3/Pin (2) = \mathbb RP^2$. 
\end{ex}

In Proposition~\ref{Wilking-Ricci} we will prove this example has
positive Ricci curvature.   In Proposition~\ref{Wilking-CovSpec},
we will prove the infinite rescaled covering spectrum is nontrivial.
These are the only facts needed to apply this example earlier in 
the paper.

Before studying the properties of this example, we explain the 
construction in more detail.   Here ${\Bbb{S}}^3 \times \Bbb{R}^4$ is
the standard isometric product with nonnegative sectional
curvature.   One takes
\be
{\Bbb{S}}^3\times \Bbb{R}^4=\{(z_1, z_2, z_3, z_4)\in \Bbb{C}^4: |z_1|^2+|z_2|^2=1\}.
\ee
We can take an ${\Bbb{S}}^1$ action on this space:
\be
e^{i\theta}(z_1, z_2, z_3, z_4)=(e^{i\theta}z_1, e^{i\theta}z_2, e^{i\theta}z_3, e^{i\theta}z_4).
\ee
Then $\tilde{M}^6= ( {\Bbb{S}}^3 \times \Bbb{R}^4 )/ {\Bbb{S}}^1$ has nonnegative
sectional curvature by O'Neill's Submersion Formula
\cite{O'Neill}.     The elements of $\tilde{M}^6$
are equivalence classes 
\be
[(z_1, z_2, z_3, z_4)]
=[(e^{i\theta}z_1,e^{i\theta}z_2,e^{i\theta}z_3,e^{i\theta}z_4)] \in \tilde{M}^6.
\ee

Recall that the Hopf sphere has ${\Bbb{S}}^2={\Bbb{S}}^3/{\Bbb{S}}^1$, where we also
view
\be
{\Bbb{S}}^3=\{(z_1, z_2)\in \Bbb{C}^2: |z_1|^2+|z_2|^2=1\}
\ee
with the same kind of circle action 
\be
e^{i\theta}(z_1, z_2)=(e^{i\theta}z_1, e^{i\theta}z_2).
\ee
The equivalence classes
\be
[(x_1+ix_3, x_2+ix_4)]=[(z_1, z_2)]
=[(e^{i\theta}z_1,e^{i\theta}z_2)] \in {\Bbb{S}}^3/{\Bbb{S}}^1
\ee
are identified with points
\be
\big(2(x_1x_2+x_3x_4), 2(x_1x_4-x_2x_3),x_1^2+x_3^2-x_2^2-x_4^2\big)
\in \{(a,b,c): a^2+b^2+c^2=1\} = {\Bbb{S}}^2\subset \R^3.
\ee
We can see that $\tilde{M}^6$ is an $\R^4$ bundle over ${\Bbb{S}}^2$
with
\be
\pi\big([(e^{i\theta}z_1, e^{i\theta}z_2),e^{i\theta}z_3, e^{i\theta}z_4)]\big)
=[(e^{i\theta}z_1, e^{i\theta}z_2)]
\ee
So it is simply connected.

We define $M^6=\tilde{M}^6/\Bbb{Z}_2$ taking the antipodal map
on the ${\Bbb{S}}^2$ base and the fibers.  First note that the antipodal
map, $g_0$, on ${\Bbb{S}}^2$:
\begin{eqnarray}
g_0[(e^{i\theta}z_1, e^{i\theta}z_2)]&=&g_0[(x_1+ix_3, x_2+ix_4)]  \\
&=&g_0\big(2(x_1x_2+x_3x_4), 2(x_1x_4-x_2x_3),x_1^2+x_3^2-x_2^2-x_4^2\big)\\
&=&\big(-2(x_1x_2+x_3x_4), -2(x_1x_4-x_2x_3), -x_1^2-x_3^2+x_2^2+x_4^2\big)\\
&=&[(x_2- i x_4,  -x_1 + x_3 i)]=[(\bar{z}_2, -\bar{z}_1)]\\
&=&[(e^{i\theta}\bar{z}_2, e^{i(\pi+\theta)}\bar{z}_1)]
\end{eqnarray}
We define the $\Bbb{Z}_2$ action $g$ on $\tilde{M}^6$ such that
$\pi\circ g=g\circ g_0$ as follows:
\be \label{M6Z2}
g[(e^{i\theta}z_1, e^{i\theta}z_2,e^{i\theta}z_3, e^{i\theta}z_4)]
=[(e^{i\theta}\bar{z}_2, e^{i(\pi+\theta)}\bar{z}_1,
 e^{i\theta}\bar{z}_4, e^{i(\pi+\theta)}\bar{z}_3)].
\ee

Since the soul of $\tilde{M}^6 = {\Bbb{S}}^3 \times \Bbb{R}^4/ {\Bbb{S}}^1$, is 
${\Bbb{S}}^3/{\Bbb{S}}^1= {\Bbb{S}}^2$, we know the soul of $M^6={\Bbb{S}}^3\times\Bbb{R}^4/\textrm{Pin}(2)=\tilde{M}^6/\Bbb{Z}_2$ is $\R P^2={\Bbb{S}}^2/\Bbb{Z}_2$.
Wilking made an argument using this soul and
the tangent cone at infinity to show the universal covering space
does not contain a line.   However, we see here that it does not
contain a line because in fact it has positive Ricci curvature:

\begin{prop}\label{Wilking-Ricci}
The Wilking Example, $M^6$, has positive Ricci curvature.
\end{prop}

\Pf
We will examine $\tilde{M}^6$.
The soul directions have positive Ricci curvature
because the soul has positive sectional
curvature and the total space has nonnegative sectional curvature. So we only need to check the fibre directions. The fibres are totally
geodesic with quotient metric $\mathbb R^4\times {\Bbb{S}}^1/{\Bbb{S}}^1$. We can
write the metric of $\mathbb R^4 \times {\Bbb{S}}^1$ as 
\be
 dr^2 +
r^2((\sigma^1)^2 + (\sigma^2)^2 + (\sigma^3)^2 + d\theta^2,
\ee
where $(\sigma^1)^2 + (\sigma^2)^2 + (\sigma^3)^2$ is the standard
left invariant metric on ${\Bbb{S}}^3$, with dual orthonormal vector fields
$X_1, X_2, X_3$ satisfy the bracket relation $[X_i, X_{i+1}] =
2X_{i+2}$ (indices are mod 3). Let $X_1$ be tangent to the Hopf
fibre direction. Then the orbit direction of ${\Bbb{S}}^1$ acts on $\mathbb
R^4\times {\Bbb{S}}^1$ is $X_1 + \frac{\partial}{\partial \theta}$.
Therefore the orthonormal basis of $\mathbb R^4\times {\Bbb{S}}^1/{\Bbb{S}}^1$ is
$\frac{\partial}{\partial r}, \bar{X}_1 = \frac{1}{\sqrt{1+r^2}}
(\frac 1r X_1 - r \frac{\partial}{\partial \theta}), \bar{X}_2 =
\frac 1r X_2, \bar{X}_3 = \frac 1r X_3$. Now we can use O'Neill's
formula to compute the sectional curvatures.
One can find
$K(\frac{\partial}{\partial r}, \bar{X}_1) = \|
[\frac{\partial}{\partial r}, \bar{X}_1]^v \|^2 =
\frac{4}{(1+r^2)^2}$, and $K(\bar{X}_2, \bar{X}_3) = \|[\bar{X}_2,
\bar{X}_3]^v\|^2 = \frac{4}{1+r^2}$. Therefore the Ricci curvature
of all directions are positive.
\qed

Intuitively, one can see that $L(g,p)\to \infty$ and $p\to \infty$
in $M^6$.  However, we will now prove it grows linearly
to obtain the following proposition:

\begin{prop}\label{Wilking-CovSpec}
The Wilking Example, $M^6$, has nontrivial infinite rescaled
covering spectrum.  
\end{prop}

\Pf
The fundamental group of $M^6$ is $\Bbb{Z}_2$.   Let $g$
be the nontrivial element, described by its action on $\tilde{M}$
in (\ref{M6Z2}).   We need only show:
\be
L_{rs}^\infty(g)=
\lim_{R\to\infty}\inf_{\,x\in M^6\setminus B_{x_0}(R)\,} \frac{d_{\tilde{M}^6}(g\tilde{x},\tilde{x})}{\,d_{M^6}(x,x_0)\,}>0.
\ee
First we compute:
\begin{eqnarray*}
d_{\tilde{M}^6}\big(g\tilde{x}, \tilde{x} \big)&=&
d_{\tilde{M}^6}\big(g[(z_1,z_2,z_3,z_4)],[(z_1,z_2,z_3,z_4)] \big)\\
&=&d_{\tilde{M}^6}\big([(e^{i\theta}\bar{z}_2, e^{i(\pi+\theta)}\bar{z}_1,
 e^{i\theta}\bar{z}_4, e^{i(\pi+\theta)}\bar{z}_3)],[(e^{i\theta}z_1, e^{i\theta}z_2,
 e^{i\theta}z_3, e^{i\theta}z_4)] \big)\\
&=&\inf_{\theta, \phi}
d_{{\Bbb{S}}^3\times \Bbb{R}^4}\big((e^{i\theta}\bar{z}_2, e^{i(\pi+\theta)}\bar{z}_1,
 e^{i\theta}\bar{z}_4, e^{i(\pi+\theta)}\bar{z}_3), (e^{i\phi}z_1, e^{i\phi}z_2, e^{i\phi}z_3, e^{i\phi}z_4)\big)\\
&=&\inf_{\theta, \phi}
d_{{\Bbb{S}}^3\times \Bbb{R}^4}\big((e^{i(\theta-\phi)}\bar{z}_2, e^{i(\pi+\theta-\phi)}\bar{z}_1,
 e^{i(\theta-\phi)}\bar{z}_4, e^{i(\pi+\theta-\phi)}\bar{z}_3), (z_1, z_2, z_3, z_4)\big)\\
&=&\inf_{\theta'}
d_{{\Bbb{S}}^3\times \Bbb{R}^4}\big((e^{i\theta'}\bar{z}_2, e^{i(\pi+\theta')}\bar{z}_1,
 e^{i\theta'}\bar{z}_4, e^{i(\pi+\theta')}\bar{z}_3), (z_1, z_2, z_3, z_4)\big)\\
&=&\inf_{\theta'}
\sqrt{\Big(d_{{\Bbb{S}}^3}\big((e^{i\theta'}\bar{z}_2, e^{i(\pi+\theta')}\bar{z}_1),(z_1, z_2)\big)\Big)^2 
+ \Big(d_{\Bbb{R}^4}\big((e^{i\theta'}\bar{z}_4, e^{i(\pi+\theta')}\bar{z}_3), 
(z_3, z_4)\big)\Big)^2}\\
&\ge&\inf_{\theta'}
d_{\Bbb{R}^4}\big((e^{i\theta'}\bar{z}_4, e^{i(\pi+\theta')}\bar{z}_3), 
(z_3, z_4)\big) \\
&=&\inf_{\theta'}\sqrt{
\big|e^{i\theta'}\bar{z}_4 -z_3\big|^2 + \big|e^{i(\pi+\theta')}\bar{z}_3 - z_4\big|^2 } \\
&\ge&\inf_{\theta'}\frac{\sqrt{2}}{2} \Big(
\big|e^{i\theta'}\bar{z}_4 -z_3\big| + \big|e^{i(\pi+\theta')}\bar{z}_3 - z_4\big| \Big) \\
&=&\inf_{\theta'}\frac{\sqrt{2}}{2} \Big(
\big|e^{-i\theta'}z_4 -\bar{z}_3\big| + \big|e^{i(\pi+\theta')}\bar{z}_3 - z_4\big| \Big) \\
&=&\inf_{\theta'}\frac{\sqrt{2}}{4} \Big(
\big|z_4 -e^{i\theta'}\bar{z}_3\big| + \big|e^{i(\pi+\theta')}\bar{z}_3 - z_4\big|
+\big|e^{-i\theta'}z_4 - \bar{z}_3\big| + 
\big|\bar{z}_3 - e^{-i(\pi+\theta')}z_4\big|
\Big)  \\
&\ge&\inf_{\theta'}\frac{\sqrt{2}}{4} \Big(
\big|e^{i\theta'}\bar{z}_3 - e^{i(\pi+\theta')}\bar{z}_3 \big|
+\big|e^{-i\theta'}z_4 - e^{-i(\pi+\theta')}z_4\big|
\Big)  \\
&=&\inf_{\theta'}\frac{\sqrt{2}}{4} \Big(
\big|\bar{z}_3 - e^{i\pi}\bar{z}_3 \big|
+\big|z_4 - e^{-i\pi}z_4\big|
\Big)  \\
&=&\frac{\sqrt{2}}{4} \Big(
\big|1 - e^{i\pi}\big| |\bar{z}_3|
+\big|1-e^{-i\pi}\big| |z_4|
\Big)  \\
&=&\frac{\sqrt{2}}{4} \Big(2 |\bar{z}_3|+2 |z_4|\Big) \\
&\ge& 
\frac{\sqrt{2}}{2} \sqrt{ |z_3|^2+ |z_4|^2}.
\end{eqnarray*}
On the other hand
\be
d_{M^6}(x,x_0) 
= \min\{d_{\tilde{M}^6}(\tilde{x}, \tilde{x}_0),
                                    d_{\tilde{M}^6}(\tilde{x}, g\tilde{x}_0) \}
                                           \ee
So by the triangle inequality
\be
\Big| d_{M^6}(x,x_0) -d_{\tilde{M}^6}(\tilde{x}, \tilde{x}_0)\Big|
\le     d_{\tilde{M}^6}(\tilde{x}_0, g\tilde{x}_0) \\                                 
\ee
Choosing $x_0$ in the soul, $\tilde{x}_0=[(w_1,w_2, 0,0)]$, so
\begin{eqnarray}
d_{\tilde{M}^6}(\tilde{x}, \tilde{x}_0)
&=&\inf_{\theta, \phi}d_{{\Bbb{S}}^3\times \Bbb{R}^4}\big((e^{i\theta}w_1, e^{i\theta}w_2,
 0,0), (e^{i\phi}z_1, e^{i\phi}z_2, e^{i\phi}z_3, e^{i\phi}z_4)\big)\\
 &\le&d_{{\Bbb{S}}^3\times \Bbb{R}^4}\big((w_1, w_2,
 0,0), (z_1, z_2, z_3, z_4)\big)\\
&\le&\sqrt{\pi^2 + |z_3|^2+|z_4|^2}.
\end{eqnarray} 

If $d_{M^6}(x, x_0)\to \infty$, then 
$d_{\tilde{M}^6}(\tilde{x}, \tilde{x}_0)\to \infty$ and
$|z_3|^2 + |z_4|^2 \to \infty$.
Thus
\begin{eqnarray}
L_{rs}^\infty(g)&=&
\lim_{R\to\infty}\inf_{\,x\in M^6\setminus B_{x_0}(R)\,} \frac{d_{\tilde{M}^6}(g\tilde{x},\tilde{x})}{\,d_X(x,x_0)\,}  \\
&\ge& 
\liminf_{|z_3|^2+|z_4|^2\to\infty} 
\frac{ \frac{\sqrt{2}}{2} \sqrt{ |z_3|^2+ |z_4|^2} }
{\sqrt{\pi^2 + |z_3|^2+|z_4|^2\,} + d_{\tilde{M}^6}(\tilde{x}_0, g\tilde{x}_0) }
=\sqrt{2}/2.
\end{eqnarray}
Thus the rescaled covering spectrum is nontrivial.
\qed


Christina Sormani

 Department of Mathematics,

CUNY Graduate Center and Lehman College

sormanic@member.ams.org

\vspace{.5cm}

Guofang Wei,

 Department of Mathematics,

 University of California Santa Barbara

 wei@math.ucsb.edu


\begin{thebibliography}{SoWei2}
 \bibitem{Anderson}
Michael~T. Anderson.
{\em On the topology of complete manifolds of nonnegative {R}icci
  curvature.}
\newblock {Topology}, 29(1):41--55, 1990.

\bibitem{Berard-Bergery} L.  B\'erard-Bergery, 
Quelques exemples de var\'et\'es riemanniennes compl\`etes non compactes \`a courbure de Ricci positive. (French. English summary) [Some examples of noncompact complete Riemannian manifolds with positive Ricci curvature]
C. R. Acad. Sci. Paris S\'er. I Math. 302 (1986), no. 4, 159--161.


\bibitem{BO} R. L. Bishop and B. O'Neill, {\em Manifolds of negative curvature}, Trans. of. Amer. Math. Soc. 
145(1969) 1--49. 

 \bibitem{BBI} D. Burago, Y. Burago, S. Ivanov, A Course
 in Metric Geometry. Graduate Studies in Mathematics Vol. 33,
 AMS, 2001.

  \bibitem{ChCo1} J. Cheeger, T. Colding,
 {\em On the structure of spaces with Ricci curvature bounded
 below I}, J. Diff. Geom. 46 (1997) 406-480.

 \bibitem{ChCo3} J. Cheeger, T. Colding, {\em On the
 structure of spaces with Ricci curvature bounded below III},
  J. Differential Geom. 54 (2000), no. 1,
37--74.

\bibitem{ChGr1971}
Jeff Cheeger and Detlef Gromoll.
{\em The splitting theorem for manifolds of nonnegative {R}icci curvature}.
\newblock {J. Differential Geometry} 6 (1971/72),  119--128.

\bibitem{Cheeger-Gromoll} J. Cheeger, D. Gromoll, {\em On the
structure of complete manifolds of nonnegative curvature},  Ann. of Math.
(2)  96  (1972), 413--443.

\bibitem{CoNb} Colding, Tobias Holck; Naber, Aaron
{\em Sharp Holder continuity of tangent cones for spaces with a lower Ricci curvature bound and applications.} 
Ann. of Math. (2) 176 (2012), no. 2, 1173--1229. 

\bibitem{CdV} Y. Colin de Verdiere, {\em Spectre du laplacien et 
 longueurs des g\'eod\'esiques p\'eriodiques. I, II.} (French) Compositio 
 Math. 27 (1973), 83--106; ibid. 27 (1973), 159--184. 

\bibitem{CCJPPWW}Jim Conant, Victoria Curnutte, Corey Jones, Conrad Plaut, Kristen Pueschel, Maria Walpole, Jay Wilkins
{\em Discrete Homotopy Theory and Critical Values of Metric Spaces},
to appear in Fundamenta Mathenaticae,  arXiv:1205.2925, May 2012.

\bibitem{deSmit-Gornet-Sutton}  de Smit, Bart; Gornet, Ruth; Sutton, Craig J. {\em Sunada's method and the covering spectrum.}  J. Differential Geom.  86  (2010),  no. 3, 501--537. 

\bibitem{deSmit-Gornet-Sutton-2}
de Smit, Bart; Gornet, Ruth; Sutton, Craig J. {\em Isospectral surfaces with distinct covering spectra via Cayley graphs. } Geom. Dedicata  158  (2012), 343--352.

\bibitem{Ding} Yu Ding, 
{\em Heat kernels and Green's functions on limit spaces.}
Comm. Anal. Geom. 10 (2002), no. 3, 475--514. 

\bibitem{DuGu}
Duistermaat, J. J.; Guillemin, V. W. 
{\em The spectrum of positive elliptic operators and periodic 
bicharacteristics.}
  Invent. Math. 29 (1975), no. 1, 39--79




\bibitem{EnWei} J. Ennis, G. Wei, 
{\em Describing the universal cover of a compact limit.}
Differential Geom. Appl. 24 (2006), no. 5, 554--562. 

 \bibitem{Gr}M. Gromov, Metric structures
 for Riemannian and non-Riemannian spaces, PM 152, Birkhauser, 1999.
 
 \bibitem{Honda} S. Honda,
{\em Bishop-Gromov type inequality on Ricci limit spaces}.
J. Math. Soc. Japan 63 (2011), no. 2, 419--442. 
 
 \bibitem{Li1986}
Peter Li.
{\em Large time behavior of the heat equation on complete manifolds with
  nonnegative {R}icci curvature.}
\newblock {Ann. of Math. (2)}, 124(1):1--21, 1986.

\bibitem{Menguy} X. Menguy, 
{\em Examples with bounded diameter growth and infinite topological type.} 
Duke Math. J. 102 (2000), no. 3, 403--412. 

 \bibitem{Mi} J. Milnor, {\em A note on curvature and
 fundamental group}, J. Diff. Geom. 2 (1968) 1-7.

\bibitem{Munn} M. Munn, 
{\em Volume growth and the topology of pointed Gromov-Hausdorff limits}. 
Differential Geom. Appl. 28 (2010), no. 5, 532--542. 

\bibitem{Myers} S.B. Myers, 
{\em Riemannian manifolds with positive mean curvature.} 
Duke Math. J. 8, (1941). 401--404. 

\bibitem{Nab} P. Nabonnand,  
Sur les var\'et\'es riemanniennes compl\`etes  \`a  courbure de Ricci positive. (French. English summary)
C. R. Acad. Sci. Paris S\'er. A-B 291 (1980), no. 10, A591 -- A593.

\bibitem{Ohta}
Shin-ichi Ohta, 
{\em On the measure contraction property of metric measure spaces.}
Comment. Math. Helv. 82 (2007), no. 4, 805--828. 

\bibitem{O'Neill}
B. O'Neill, {\em The fundamental equations of a submersion}, Michigan Math. J. 13 (1966), 459--469. 

\bibitem{Perelman1} G. Perelman, {\em Proof of the soul conjecture of Cheeger and Gromoll}, J. Differential Geom. 40 (1994), no. 1, 209--212. 

\bibitem{Perelman2} G. Perelman, {\em Alexandrov's spaces with curvatures bounded from below II}, perpetual preprint.

\bibitem{Plaut-Wilkins:1205:1137} C. Plaut and J. Wilkins, {\em
Discrete homotopies and the fundamental group},  Advances in Mathematics 232 (2013), 271--294.

\bibitem{Sharafutdinov} V. A. Sharafutdinov, {\em The Pogorelov-Klingenberg theorem for manifolds homeomorphic to $\mathbb R^n$}, Sib. Math. Zh. 18 (1977) 915-925.

 \bibitem{Sormani-group}
Christina Sormani.
\newblock Nonnegative {R}icci curvature, small linear diameter growth and
  finite generation of fundamental groups.
\newblock {\em J. Differential Geom.}, 54(3):547--559, 2000.

 \bibitem{So-loops} C. Sormani, {\em On Loops Representing Elements of
 the Fundamental Group of a Complete Manifold with Nonnegative Ricci
 Curvature},  Indiana Journal of Mathematics, 50 (2001) no. 4, 1867-1883.

\bibitem{Sormani-length} C. Sormani, {\em Convergence and the length spectrum.}
Advances in Mathematics 213 (2007), no. 1, 405--439. 

 \bibitem{SoWei1} C. Sormani and G. Wei, {\em Hausdorff Convergence and
 Universal Covers}, Transactions of the American Mathematical Society 353
 (2001) 3585-3602.

  \bibitem{SoWei2} C. Sormani and G. Wei, {\em Universal Covers for
 Hausdorff Limits of Noncompact Spaces}, Transactions of the American
 Mathematical Society 356 (2004) no. 3, 1233-1270.

\bibitem{SoWei3}  C. Sormani and G. Wei, {\em The Covering
 Spectrum of a Compact Length Space},
 Journal of Differential Geometry, Vol 67 (2004) 35-77.

\bibitem{SoWei4}  C. Sormani and G. Wei, {\em The Cut-off Covering
 Spectrum},  Tran. AMS. 362, no. 5 (2010), 2339-2391. 

\bibitem{Sp} E. Spanier, Algebraic Topology, McGraw-Hill, Inc., 1966.

\bibitem{Wei1988}
Guofang Wei.
\newblock Examples of complete manifolds of positive {R}icci curvature with
  nilpotent isometry groups.
\newblock {\em Bull. Amer. Math. Soci.}, 19(1):311--313, 1988.

\bibitem{Wilking2000}
Burkhard Wilking.
\newblock On fundamental groups of manifolds of nonnegative curvature.
\newblock {\em Differential Geom. Appl.}, 13(2):129--165, 2000.

\bibitem{Wilkins-thesis} J. Wilkins, {\em Discrete Geometric Homotopy
Theory and Critical Values of Metric Spaces}, Doctoral Dissertation,
University of Tennesee Knoxville, 2011.

\bibitem{Wilkins-1208.3494} J Wilkins, {\em The Revised and Uniform Fundamental Groups and Universal Covers of Geodesic Spaces}.
Topology and its Applications 160 (2013), no. 6, 812--835.

\bibitem{XWY} S. Xu, Z. Wang, F. Yang,
{\em On the fundamental group of open manifolds with nonnegative Ricci curvature.} Chinese Ann. Math. Ser. B 24 (2003), no. 4, 469--474. 
 

 \end{thebibliography}
 \end{document}